\documentclass[12pt]{article}
\usepackage{mathtools}
\usepackage[utf8]{inputenc}
\usepackage{amssymb,amsmath,amsfonts,eucal,mathrsfs,amsthm}
\usepackage[colorinlistoftodos]{todonotes}
\usepackage{xcolor}
\usepackage{enumitem}
\newtheorem{theorem}{Theorem}
\newtheorem{proposition}[theorem]{Proposition}
\newtheorem{lemma}[theorem]{Lemma}

\newtheorem{corollary}[theorem]{Corollary}

\theoremstyle{definition}
\newcommand{\R}{\mathbb{R}}
\newcommand{\Z}{\mathbb{Z}}
\newcommand{\Q}{\mathbb{Q}}
\newcommand{\Sf}{\mathbb{S}}

\newcommand{\Hy}{\mathbb{H}}

\newcommand{\spa}{\mbox{span}}

\newcommand{\Ric}{\mbox{Ric}}

\newcommand{\trace}{\mbox{tr\,}}

\def\<{{\langle}}
\def\>{{\rangle}}
\def\B{\mathcal{B}}
\def\CP{\mathord{\mathbb C}\mathord{\mathbb P}}

\def\n{\nabla}

\def\a{\alpha}

\def\be{\begin{equation} }
\def\ee{\end{equation} }

\begin{document}

\title{Geometric and topological rigidity of pinched submanifolds}
\author{Theodoros Vlachos}
\date{}
\maketitle

\renewcommand{\thefootnote}{\fnsymbol{footnote}} 
\renewcommand{\thefootnote}{\arabic{footnote}} 

\renewcommand{\thefootnote}{\fnsymbol{footnote}} 
\renewcommand{\thefootnote}{\arabic{footnote}}

\begin{abstract}
We investigate the geometry and topology of compact submanifolds of 
arbitrary codimension in space forms satisfying a certain pinching 
condition involving the length of the second fundamental form and 
the mean curvature. We prove that this pinching condition either 
forces homology to vanish in a range of intermediate dimensions, 
or completely determines the submanifold up to congruence. The 
results are sharp and extend previous results due to several authors 
without imposing any further assumption on the mean curvature. 
\end{abstract}

\section{Introduction}

A fundamental problem in differential geometry is to investigate the 
relationship between geometry and topology of Riemannian manifolds. 
From the point of view of submanifold theory, it should be interesting 
to investigate how the geometry or the topology of submanifolds 
of space forms is affected by pinching 
conditions involving intrinsic and extrinsic curvature invariants. 
For minimal submanifolds of spheres with 
sufficiently pinched second fundamental form, a fundamental 
result was first given by Simons in his seminal 
paper \cite{Sim}. Afterwards, Chern, do Carmo and Kobayashi 
\cite{CdCK} proved a well known rigidity theorem. Since then, 
their work has inspired plenty of interesting pinching results, 
see for instance \cite{adc, FX,HW,Le1,OV,SX, V2, XG}. 

Lawson and Simons \cite{LS} showed that certain 
bounds on the second fundamental form for submanifolds of 
spheres imply vanishing of homology groups with integer 
coefficients. They used the second variation of area to 
rule out stable minimal currents in certain dimensions. 
Since one can minimize area in a homology class (see \cite{FF}), 
this trivializes integral homology. 

The aim of this paper is to study the geometric and 
topological rigidity of submanifolds under a sharp pinching 
condition involving the length of the second fundamental form 
 and the mean curvature. Throughout this paper, 
$S$ denotes the {\it squared length} of the second 
fundamental form $\a_f$ of an isometric immersion $f$, 
while the \textit{mean curvature} is defined as the length 
$H=\| \mathcal H \|$ of the 
\textit{mean curvature vector field} given by 
$\mathcal H=(\trace\alpha_f)/n$, where 
$\mathrm{tr}$ means taking the trace.

The choice of the pinching condition is 
inspired by the standard immersion of a torus 
$\mathbb T^n_k(r)=\mathbb S^k(r) \times \mathbb S^{n-k}(\sqrt{1-r^2}), 1\leq k\leq n-1$,
into the unit sphere $\mathbb S^{n+1}$, where $\mathbb S^k(r)$ 
denotes the $k$-dimensional sphere of radius $r<1$. The 
principal curvatures are $\sqrt{1-r^2}/r$ and $-r/\sqrt{1-r^2}$ 
of multiplicities $k$ and $n-k$, respectively. A direct computation 
gives that the squared length $S$ of the second 
fundamental form satisfies 
\begin{equation*}
S=n+\frac{n^3H^2}{2k(n-k)}\pm \frac{n|n-2k|}{2k(n-k)}H\sqrt{
n^2H^2+4k(n-k)},
\end{equation*}
where the sign $+$ or $-$ is chosen according to $r\leq \sqrt{k/n}$, or 
$r\geq \sqrt{k/n}$, respectively. 
Then, we have $S=a(n,k,H,1)$ if $r\geq \sqrt{k/n}$ or $n=2k$ and 
$S>a(n,k,H,1)$ if otherwise, where the function $a$ is given by 
\begin{align*}
a(n,k,t,c)
&=nc(2-{\rm {sgn}}(c))+\frac{n^3t^2}{2k(n-k)}\\
&-\frac{n|n-2k|}{2k(n-k)}t\sqrt{
n^2t^2+4ck(n-k)(2-{\rm {sgn}}(c))},\;t, c\in \mathbb{R},
\end{align*}
with $n^2t^2+4ck(n-k)(2-{\rm {sgn}}(c))\geq0$.

Motivated by this example, we study isometric immersions 
$f\colon M^n\to \Q_c^{n+m}$ that satisfy the pinching condition
\be\label{1}
S\leq a(n,k,H,c) \tag{$\ast$}
\ee
pointwise, where $k$ is an integer with $1\leq k\leq n/2$. Here 
$\Q_c^{n+m}$ denotes the $(n+m)$-dimensional complete simply 
connected space form of constant curvature $c$. For simplicity we 
assume that $c\in\{0,\pm1\}$, unless otherwise stated. Thus 
$\Q_c^{n+m}$ is the Euclidean space $\R^{n+m}$, the unit 
sphere $\Sf^{n+m}$ or the hyperbolic space $\Hy^{n+m}$, 
if $c=0,1,-1$, respectively.

The pinching condition \eqref{1} has been studied mainly for 
the lowest allowed $k$ (see for instance \cite{WXia,Xu,XG,ZZ}). 
Shiohama and Xu \cite{SX} proved that compact submanifolds in space 
forms of nonnegative curvature are homeomorphic to a sphere provided 
that \eqref{1} holds as strict inequality at any point for $k=1$. 
In \cite{V2}, a homology vanishing result was obtained for submanifolds in 
space forms of nonnegative curvature that satisfy strict inequality in 
\eqref{1} pointwise for an integer $1\leq k\leq n-1$. The approach therein 
was based on the nonexistence results of stable currents due to Lawson 
and Simons \cite{LS} under certain upper bounds of the second 
fundamental form.

In the present paper, we investigate the geometric and topological rigidity 
of submanifolds that satisfy the pinching condition \eqref{1}, where we 
no longer ask to be strict. We are able to show that the pinching condition 
either forces homology to vanish in a range of intermediate dimensions, or 
completely determines the pinched submanifold up to congruence. 
Throughout the paper all submanifolds under consideration are assumed 
to be connected. We recall that a submanifold is called substantial if it is not 
contained in a proper totally geodesic submanifold of the ambient space. 

\begin{theorem}\label{th1}
Let $f\colon M^n\to\Q_c^{n+m}$, 
$n\geq 5$, be a substantial isometric immersion of a compact 
oriented Riemannian manifold. Assume that the inequality \eqref{1} is 
satisfied{\footnote{Notice that if $c=-1$ and 
$H\geq\sqrt{3}$, then $n^2H^2+4ck(n-k)(2-{\rm {sgn}}(c))\geq0$.}} for an 
integer $2\leq k\leq n/2$ at any point. If $c=-1$, suppose further that the 
mean curvature satisfies $H\geq\sqrt{3}$ everywhere. Then, either the 
homology groups of $M^n$ satisfy 
$$
H_p(M^n;\Z)=0\;\,\text{for all}\;\,k\leq p\leq n-k
\;\,\text{and}\;\, H_{k-1}(M^n;\Z)=\Z^{\beta_{k-1}(M)},
$$
where $\beta_{k-1}(M)$ 
is the $(k-1)$-th Betti number, or equality holds in \eqref{1} at any point and 
one of the following assertions holds:
\item[(i)] $M^n$ is isometric to a Clifford torus 
$\mathbb T^n_p(\sqrt{p/n}),k\leq p\leq n/2$, 
and $f$ is the standard minimal embedding into $\Sf^{n+1}$.
\item[(ii)] $M^n$ is isometric 
to a torus $\mathbb T^n_k(r)$ with $r>\sqrt{k/n}$ 
and $f$ is the standard embedding into $\Sf^{n+1}$.
\item[(iii)] $M^n$ is isometric 
to a torus $\Sf^k(r)\times \Sf^k(\sqrt{R^2-r^2})$ and $f$ is a composition 
$f=j\circ g$, where $g\colon M^n\to\Sf^{n+1}(R)$ is the standard embedding 
of the torus $\Sf^k(r)\times \Sf^k(\sqrt{R^2-r^2})$ into a sphere $\Sf^{n+1}(R)$ 
and $j\colon\Sf^{n+1}(R)\to\Q_c^{n+2}$ is an umbilical inclusion with $c=0,1$ 
and $R<1$ if $c=1$.
\end{theorem}

The class of submanifolds as in part $(i)$ of Theorem \ref{th1} is nonempty. 
Indeed, umbilical submanifolds do satisfy \eqref{1} for any $2\leq k\leq n/2$ 
at any point. Besides, one can readily verify that the standard embedding of 
the torus $\mathbb T^n_{k-1}(r)$ into $\Sf^{n+1}$ satisfies \eqref{1} for 
$2\leq k\leq n/2$ at any point, if $r$ is close to $\sqrt{(k-1)/n}$ and 
$r\geq\sqrt{(k-1)/n}$. Moreover, the homology groups is as required in part 
$(i)$ of Theorem \ref{th1}, since 
$H_p(\mathbb T^n_{k-1}(r);\Z)=0\;\,\text{for all}\;\,k\leq p\leq n-k$ and 
$H_{k-1}(\mathbb T^n_{k-1}(r);\Z)=\Z$, or 
$H_{k-1}(\mathbb T^n_{k-1}(r);\Z)=\Z\oplus\Z$ if $n\neq2(k-1)$, or 
$n=2(k-1)$, respectively. Moreover, Wallach \cite{W} 
constructed a minimal embedding $\psi_k\colon\CP^k_{2k/(k+1)}\to\Sf^{k(k+2)-1}$ 
of the complex projective space of constant holomorphic curvature $2k/(k+1)$ 
which satisfies \eqref{1} only if $k=2$. This shows that Theorem \ref{th1} 
doesn't hold for $n=4$ since $H_2(\CP^2;\Z)=\Z$. 

For submanifolds satisfying the pinching condition \eqref{1} for $k=1$, 
we prove the following stronger result.

\begin{theorem}\label{th2}
Let $f\colon M^n\to\Q_c^{n+m}$, 
$n\geq 3$, be a substantial isometric immersion of a compact 
oriented Riemannian manifold. Assume that the inequality \eqref{1} 
is satisfied for $k=1$ at any point. If $c=-1$, suppose further that the 
mean curvature satisfies $H\geq\sqrt{3}$ everywhere.
Then, either $M^n$ is homeomorphic to $\Sf^n$ if $n\geq5$, diffeomorphic 
to a spherical space form if $n=3$, its universal cover is homeomorphic to 
$\Sf^4$ if $n=4$, or equality holds in \eqref{1} for $k=1$ at any point and 
one of the following assertions holds: 
\item[(i)] $M^n$ is isometric to a Clifford torus 
$\mathbb T^n_p(\sqrt{p/n}),1\leq p\leq n-1$, 
and $f$ is the standard minimal embedding into $\Sf^{n+1}$.
\item[(ii)] $M^n$ is isometric to a torus $\mathbb T^n_1(r) $ 
with $r>1/\sqrt{n}$ and $f$ is the standard embedding into $\Sf^{n+1}$. 
\item[(iii)] $M^n$ is isometric to the complex projective plane of constant 
holomorphic curvature $4/3$ and $f$ is the minimal embedding 
$\psi_2\colon \CP^2_{4/3}\to \Sf^7$ due to Wallach. 
\end{theorem}

In case where the manifold $M^n$ is topologically a sphere in the above 
result and if $n=5,6,12,56,61$, then $M^n$ is diffeomorphic to $\Sf^n$ 
for these dimensions the differentiable structure of the sphere is unique; 
see Corollary $1.15$ in \cite{WX}. The above result strengthens the 
sphere theorem due to Shiohama and Xu \cite{SX}. Moreover, it improves 
and extends results in \cite{AC,OV,WXia,Xu,XG,ZZ} for any codimension 
without imposing any further assumption either on the mean curvature or 
on the fundamental group of the submanifold.

Howard and Wei in Theorem 7 of \cite{HW} investigated the topology of submanifolds 
$f\colon M^n\to\Q_c^{n+m}, c\geq0$, satisfying the strict inequality 
$$
S<n^2H^2/(n-k)+2kc
$$ for an integer $1\leq k\leq n/2$ at any point. 
It is worth pointing out that our pinching condition \eqref{1} is weaker than 
the one of Howard and Wei, since $n^2H^2/(n-k)+2kc\leq a(n,k,H,c)$. 
Thus Theorems \ref{th1} and \ref{th2} sharpen Theorem 5.7 in \cite{HW} 
and provide a partial answer to a problem raised therein for submanifolds 
in space forms with curvature $c<0$. Stronger pinching conditions have 
been studied by means of the mean curvature flow 
(see for instance \cite{AB,LXYZ}).

The following direct consequence of Theorem \ref{th2} improves previous 
results for minimal submanifolds in spheres (see \cite{CdCK,Le1,OV}).

\begin{corollary}\label{cor}
Let $f\colon M^n\to\Sf^{n+m}$, 
$n\geq 3$, be a substantial isometric minimal immersion of a compact 
oriented Riemannian manifold. If $S\leq n$ at any point, then either 
$M^n$ is homeomorphic to $\Sf^n$ if $n\geq5$, diffeomorphic to a 
spherical space form if $n=3$, its universal cover is homeomorphic to 
$\Sf^4$ if $n=4$, or $S=n$ at any point and one of the following holds:
\item[(i)] The manifold $M^n$ is isometric to a Clifford torus 
$\mathbb T^n_p(\sqrt{p/n}),1\leq p\leq n-1$, 
and $f$ is the standard embedding into $\Sf^{n+1}$.
\item[(ii)] $M^n$ is isometric to the complex projective plane of constant 
holomorphic curvature $4/3$ and $f$ is the minimal embedding 
$\psi_2\colon \CP^2_{4/3}\to \Sf^7$ due to Wallach.
\end{corollary}

The following is a consequence of Theorem \ref{th1}.

\begin{corollary}\label{k=2}
Let $f\colon M^n\to\Q_c^{n+m}$, 
$n\geq 5$, be a substantial isometric immersion of a compact 
oriented Riemannian manifold. Assume that the inequality \eqref{1} 
is satisfied for $k=2$ at any point. If $c=-1$, suppose further that the 
mean curvature satisfies $H\geq\sqrt{3}$ everywhere. If the fundamental 
group of $M^n$ is finite, then either $M^n$ is homeomorphic to $\Sf^n$, 
or equality holds in \eqref{1} for $k=2$ at any point and one of the following holds: 
\item[(i)] $M^n$ is isometric to a Clifford torus 
$\mathbb T^n_p(\sqrt{p/n}),2\leq p\leq n-2$, 
and $f$ is the standard minimal embedding into $\Sf^{n+1}$.
\item[(ii)] $M^n$ is isometric to a torus $\mathbb T^n_2(r) $ 
with $r>\sqrt{2/n}$ and $f$ is the standard embedding into $\Sf^{n+1}$. 
\end{corollary}

\section{The pinching condition}

The result of Lawson and Simons \cite{LS} mentioned in the introduction 
was later strengthened by Elworthy and Rosenberg \cite[p.\ 71]{ER} 
without requiring the bound on the second fundamental form to be strict 
at all points of the submanifold. The case when the ambient space is a 
hyperbolic space was considered in \cite{FX}. In this section, we state 
their theorem and then give an auxiliary result that establishes 
the relation between our pinching condition \eqref{1} and the 
inequality \eqref{LS} below required in their result.

\begin{theorem}(\cite{ER,FX,LS, X})\label{er} Let $f\colon M^n\to\Q_c^{n+m}$, 
$n\geq 4$, 
be an isometric immersion of a compact Riemannian manifold and $p$ be an integer
such that $1\leq p\leq n-1$. Assume that at any point $x\in M^n$ 
and for \emph{any} orthonormal basis $\{e_1,\ldots,e_n\}$ of $T_xM$ 
the second fundamental form $\alpha_f\colon TM\times TM\to N_fM$ 
satisfies
\be\label{LS}
\sum_{i=1}^p\sum_{j=p+1}^n\big(2\|\a_f(e_i,e_j)\|^2
-\<\a_f(e_i,e_i),\a_f(e_j,e_j)\>\big)\leq cp\big(n-p\,{\rm {sgn}}(c)\big) \tag{$\ast\ast$}. 
\ee
If there is a point where the inequality \eqref{LS} is strict for any orthonormal 
basis of the tangent space at that point, then the homology groups satisfy
$H_p(M^n;\mathbb{Z})=H_{n-p}(M^n;\mathbb{Z})=0$.
\end{theorem}
\proof
For $c\geq0$ the above result was proved by Lawson and Simons \cite{LS} if 
strict inequality holds in \eqref{LS} at any point (see also \cite{X} for $c=0$). 
Elworthy and Rosenberg \cite[p.\ 71]{ER} observed that it also holds by only 
requiring the bound to be strict at some point of the submanifold. 

In the case $c<0$, the theorem was proved by Fu and Xu \cite{FX} 
under the condition that strict inequality holds at any point. 
A similar  observation of Elworthy and Rosenberg \cite[p.\ 71]{ER} 
also applies to the case $c<0$, provided there exists a point 
where the inequality \eqref{LS} is strict. 

We briefly outline the structure of the proof. Recall that 
a rectifiable $p$-current $\mathscr G$ is said to be stable if, 
for any smooth vector field $V\in\mathcal X(M)$, there exists 
$\varepsilon>0$ such that the mass $M(\mathscr G)$ satisfies 
$$
M(\mathscr G)\leq M(\mathscr \phi_{t_*}\mathscr G)\; {\text {for}}
\;|t|<\varepsilon,
$$
where $\phi_t\colon M\to M$ is the one-parameter group of 
diffeomorphisms generated by $V$, and $\mathscr \phi_{t_*}\mathscr G$ 
is the current defined as a linear functional on the space of $p$-forms of 
$M$ by 
$\mathscr \phi_{t_*}\mathscr G(\omega)=\mathscr G ({\phi_{t}}^*\omega)$. 

Viewing the hyperbolic space $\mathbb H _c^{n+m}$ as a 
submanifold of the Lorentz-Minkowski space 
$\mathbb L_1^{n+m+1}$ via the standard inclusion, we 
consider the space $\mathcal V$ consisting of gradients of functions 
$f_v=\<v,f\>$, where $v$ is any vector in $\mathbb L_1^{n+m+1}$. 

To each $p$-current $\mathscr G$, we associate a quadratic form 
$Q_{\mathscr G}$ on $\mathcal V$ as follows. For each 
$V\in\mathcal V$, we define  
$$
Q_{\mathscr G}(V)=\frac{d^2}{dt^2}M({\phi_{t}}_*\mathscr G)\vert_{t=0},
$$
where $\phi_t$ is the flow generated by $V$. There is a natural 
isomorphism between $\mathcal V$ and $\mathbb L_1^{n+m+1}$. 
This isomorphism induces a natural inner product on $\mathcal V$ 
and $Q_{\mathscr G}$ can be regarded as a quadratic form on 
$\mathbb L_1^{n+m+1}$.

Lawson and Simons computed the first and second 
variation formulas for the mass in Theorem 1 of \cite{LS}.  
Based on this, and following the computation in 
(see \cite[pp. 815-816]{FX}), we obtain 
$$
\trace Q_{\mathscr G}=
\int\big(\mathbf{B}(\vec{\mathscr G}\,)-p(n+p)c\big)\,d\,\|\mathscr G\|.
$$
Here $\mathbf{B}$ is a function on the 
Grassmannian $G_{p,n-p}$ of $p$-planes in $T_xM$ 
at each point $x\in M^n$, defined by
$$
\mathbf{B}(\xi)=\sum_{i=1}^p\sum_{j=p+1}^n\big(2\|\a_f(e_i,e_j)\|^2
-\<\a_f(e_i,e_i),\a_f(e_j,e_j)\>\big), 
$$
where $\xi=e_1\wedge\dots\wedge e_p$ and
$\{e_1,\dots,e_p, e_{p+1},\dots,e_n\}$ is an orthonormal basis of 
$T_xM$. 
 
Clearly, if there is a point where the inequality \eqref{LS} is strict 
for every orthonormal basis of the tangent space at that point, then
we obtain $\trace Q_{\mathscr G}<0$. This rules out 
the existence of stable currents in dimensions $p$ and $n-p$. The 
vanishing of the homology groups in these dimensions then follows 
from the existences theorems of Federer and Fleming \cite{FF}, which 
assert that any nontrivial integral homology class corresponds to a 
stable current. \qed

\vspace{1ex}

Let $f\colon M^n\to\Q_c^{n+m}$ be an isometric immersion. 
We recall that a vector $\eta$ in the normal space $N_fM(x)$ 
is called a \emph{principal normal} of the isometric immersion 
$f$ at a point $x\in M^n$ if the tangent subspace
$$
E_\eta(x)=\left\{X\in T_xM:\alpha_f(X,Y)
=\<X,Y\>\eta\;\,\text{for all}\;\,Y\in T_xM\right\}
$$
is nontrivial. The dimension of this subspace is called the 
\emph{multiplicity} of $\eta$. If the multiplicity is at least two, then 
$\eta$ is called a \emph{Dupin principal normal}. 
The \emph{relative nullity} subspace $\mathcal D_f(x)$ of $f$ at a point 
$x\in M^n$ is the kernel of the second fundamental form at this point, 
namely
$$
\mathcal D_f(x)=\left\{X\in T_xM:\alpha_f(X,Y)
=0\;\,\text{for all}\;\,Y\in T_xM\right\}.
$$
The \emph{first normal space} of the isometric immersion $f$ at 
$x\in M^n$ is defined as 
$$
N_1f(x)={\rm {span}}\left\{\a_f(X,Y):X,Y\in T_xM\right\}.
$$
Observe that the orthogonal 
complement of $N_1f(x)$ in $N_fM(x)$ is 
$$
N_1^\perp f(x)=\{\xi\in N_fM(x):A_\xi=0\},
$$
where $A_{\xi}$ is the \emph{shape operator} of $f$ associated 
to any normal vector $\xi$. The points where $f$ is \emph{umbilical} 
are the zeros of the traceless part $\Phi=\a_f-\<\cdot,\cdot\>\mathcal H$ 
of the second fundamental form. 

Next we analyse the relation between the inequalities \eqref{1}
and \eqref{LS}. 

\begin{proposition}\label{propu}
Let $f\colon M^n\to\Q_c^{n+m},n\geq 4$, be an 
isometric immersion such that the inequality \eqref{1} is satisfied at a 
point $x\in M^n$ for an integer $1\leq k\leq n/2$. If $c=-1$, suppose 
further that $H(x)\geq\sqrt{3}$. Then, the inequality \eqref{LS} is 
satisfied for any integer $k\leq p\leq n/2$ with $p>1$ and is strict 
if the inequality \eqref{1} is strict at $x$. 

Suppose now that equality holds in \eqref{LS} for a certain orthonormal 
basis $\{e_i\}_{1\leq i\leq n}$ of $T_xM$ and an integer $k\leq p\leq n/2$ 
with $p>1$. Then, the following assertions hold at $x$:
\item[(i)] If $n\geq5$, then $f$ has flat normal bundle at $x$ and the 
following facts hold:
\begin{itemize}
\item[(i1)] If $f$ is not umbilical at $x$, then there exist two distinct Dupin 
principal normals $\eta_1,\eta_2$ at $x$ such that $
E_{\eta_1}(x)=\spa\left\{e_1,\dots,e_p\right\}$ and 
$E_{\eta_2}(x)=\spa\left\{e_{p+1},\dots,e_n\right\}$. 
\item[(i2)] If $f$ is umbilical at $x$, then either $c=0$ and $f$ is totally 
geodesic at $x$, or $c=-1$ and $H(x)=\sqrt{3}$. 

Moreover, if $k<n/2$, then we have:
\item[(i3)] If $c=-1$, then $n=2p$ and $f$ is umbilical at $x$ with 
$H(x)=\sqrt{3}$.
\item[(i4)] If $c\geq0, H(x)\neq0$ and $f$ is not umbilical at $x$, 
then $p=k$ and $N_1f(x)={\rm {span}}\left\{\mathcal H(x)\right\}$.
\item[(i5)] If $c\geq0, H(x)=0$ and $f$ is not umbilical at $x$, then 
$c=1$ and there exists a unit vector $\xi$ such that 
$N_1f(x)={\rm {span}}\left\{\xi\right\}$. Moreover, the eigenvalues of the 
shape operator $A_\xi$ at $x$ are $\sqrt{(n-p)/p}$ and 
$-\sqrt{p/(n-p)}$ with corresponding multiplicities $p$ and $n-p$.
\end{itemize}
\item[(ii)]If $n=4$ and $p=2$, then there are normal vectors $\eta_j$, 
$j=1,2$, such that 
\be\label{Vj}
\pi_j\circ A_\xi|_{V_j}=\<\xi,\eta_j\>Id\;\,
\text{for any}\;\,\xi\in N_fM(x),\nonumber
\ee
where $Id$ is the identity map on the tangent space at $x$, 
$V_1=\spa\{e_1,e_2\},V_2=\spa\{e_3,e_4\}$ and $\pi_j$ denotes 
the projection onto $V_j,j=1,2$. 
\end{proposition}

The following lemma was proved in \cite{V2}.

\begin{lemma}\label{lem}
Let $x_1,\dots,x_n$ be real numbers and $p,q$ 
be positive integers such that $p+q=n$. We set 
$$
\sigma
_1=\sum_{r=1}^nx_r \;\;\text{and}\;\; \sigma _2=\sum_{r=1}^nx_r^2.
$$
Then, we have
\begin{equation*}
-\sum_{i=1}^px_i\sum_{j=p+1}^nx_j\leq \frac{pq}n\Big(\sigma _2-
\frac{2}{n}\sigma _1^2+\frac{|p-q|}{n\sqrt{pq}}\left\vert \sigma
_1\right\vert \sqrt{n\sigma _2-\sigma _1^2}\Big).
\end{equation*}
Moreover, if equality holds in the above inequality, then 
$x_1=\dots=x_p$ and $x_{p+1}=\dots=x_n$.
\end{lemma}

We need the following elementary fact.

\begin{lemma}\label{lemn}
Let $f\colon M^n\to\Q_c^{n+m},n\geq 4$, be an 
isometric immersion. If $c=-1$, we suppose that $H(x)\geq\sqrt{3}$ 
at a point $x\in M^n$. Then, the inequality \eqref{1} holds at $x$ 
for an integer $1\leq k\leq n/2$ if and only if $\|\Phi\|\leq r(n,k,H,c)$ 
at $x$, where $r(n,k,H,c)$ is the nonnegative root of the polynomial 
\be\label{pol}
P_{n,k,H,c}(t)=t^2+\frac{n(n-2k)}{\sqrt{nk(n-k)}}Ht-n\big(H^2+c\big(2-{\rm {sgn}}(c)\big)\big).
\ee
Moreover, equality holds in \eqref{1} if and only if $\|\Phi\|= r(n,k,H,c)$.
\end{lemma}
\proof
Observe that the assumption on the mean curvature implies that 
$$
n^2H^2+4ck(n-k)\big(2-{\rm {sgn}}(c)\big)\geq0.
$$
Hence, the polynomial given by \eqref{pol} has real roots, one of 
them being nonnegative. Since $\|\Phi\|^2=S-nH^2$, the pinching 
condition \eqref{1} is equivalent to the inequality 
$$
\|\Phi\|^2\leq a(n,k,H,c)-nH^2.
$$
A direct computation shows that 
$$
a(n,k,H,c)=(r(n,k,H,c))^2+nH^2,
$$
where $r(n,k,H,c)$ is given by 
\begin{equation*}
r(n,k,H,c)=\frac{\sqrt{n}}{2\sqrt{k(n-k)}}\left|\sqrt{n^2H^2+4ck(n-k)\big(2-{\rm {sgn}}(c)\big)}-(n-2k)H\right|.
\end{equation*}
It is clear that $r(n,k,H,c)$ is the nonnegative root of the polynomial 
given by \eqref{pol}.\qed
\vspace{1.5ex}

\noindent\emph{Proof of Proposition \ref{propu}:}
Let $\{e_i\}_{1\leq i\leq n}$ be an arbitrary orthonormal basis 
of $T_xM$. We choose an orthonormal basis 
$\{\xi_\a\}_{1\leq\a\leq m}$ of $N_fM(x)$ such that $\mathcal{H}(x)=H(x)\xi_1$. 
Let $A_\a,1\leq\a\leq m$, be the corresponding shape operators. 
Setting 
$h_{ij}^{(\a)}=\<A_{\a}e_i,e_j\>, 1\leq\a\leq m$, we have 
$$
\sum\limits_{\ell=1}^nh_{\ell\ell}^{(1)}=nH \;\;{\text {and}}\;\;
\sum\limits_{\ell=1}^nh_{\ell\ell}^{(\alpha)}=0\ \ \text{for all}\ \ 2\leq \alpha\leq m.
$$
Let $p$ be an integer such that $k\leq p\leq n/2$ with $p>1$. From 
Lemma \ref{lem}, it follows that 
\begin{equation}\label{in1u}
-\sum\limits_{i=1}^ph_{ii}^{(\alpha)}\sum\limits_{j=p+1}^nh_{jj}^{(\alpha)}\leq \frac{p(n-p)}{n}\sum\limits_{\ell=1}^{n}(h_{\ell\ell}^{(\alpha)})^2\ \ \text{for all}\ \ 
2\leq \alpha\leq m,
\end{equation}
and 
\begin{align}\label{in2u}
&-\sum\limits_{i=1}^ph_{ii}^{(1)}\sum\limits_{j=p+1}^nh_{jj}^{(1)}\nonumber\\
&\leq \frac{p(n-p)}{n}\Big(\sum\limits_{\ell=1}^{n}(h_{\ell\ell}^{(1)})^2
-2nH^2
+\frac{n|n-2p|H}{\sqrt{np(n-p)}}\sqrt{\sum\limits_{\ell=1}^n(h_{\ell\ell}^{(1)})^2-nH^2}\Big).
\end{align}
Using \eqref{in1u}, \eqref{in2u} and since 
\begin{align*}
&\sum_{i=1}^p\sum_{j=p+1}^n\big(2\|\a_f(e_i,e_j)\|^2
-\<\a_f(e_i,e_i),\a_f(e_j,e_j)\>\big) \\
&=\sum\limits_{\alpha=1}^m\Big(
2\sum\limits_{i=1}^{p}\sum\limits_{j=p+1}^n(h_{ij}^{(\alpha)})^2-\sum\limits_{i=1}^ph_{ii}^{(\alpha)}\sum\limits_{j=p+1}^nh_{jj}^{(\alpha)}\Big),
\end{align*}
we obtain 
\begin{align}\label{in3u}
&\sum_{i=1}^p\sum_{j=p+1}^n\big(2\|\a_f(e_i,e_j)\|^2
-\<\a_f(e_i,e_i),\a_f(e_j,e_j)\>\big)\nonumber\\
&\leq 2\sum\limits_{i=1}^p\sum\limits_{j=p+1}^n\sum\limits_{\alpha
=1}^m(h_{ij}^{(\alpha)})^2+\frac{p(n-p)}{n}\sum\limits_{\alpha
=2}^m\sum\limits_{\ell=1}^n(h_{\ell\ell}^{(\alpha)})^2\nonumber\\
&+\frac{p(n-p)}{n}\Big(\sum\limits_{\ell=1}^n(h_{\ell\ell}^{(1)})^2-2nH^2+\frac{n|n-2p|H}{\sqrt{np(n-p)}}\sqrt{
\sum\limits_{\ell=1}^n(h_{\ell\ell}^{(1)})^2-nH^2}\Big).
\end{align}
Observe that if $n>4$ then $p(n-p)>n$ for any $2\leq p\leq n-2$. 
Then, the inequality \eqref{in3u} yields 
\begin{align}\label{in4u}
\sum_{i=1}^p\sum_{j=p+1}^n&\big(2\|\a_f(e_i,e_j)\|^2
-\<\a_f(e_i,e_i),\a_f(e_j,e_j)\>\big)\nonumber\\
&\leq\frac{p(n-p)}{n}\Big(2\sum\limits_{i=1}^p\sum\limits_{j=p+1}^n\sum\limits_{\alpha
=1}^m(h_{ij}^{(\alpha)})^2+\sum\limits_{\alpha=1}^m\sum\limits_{\ell=1}^n(h_{\ell\ell}^{(\alpha)})^2 \nonumber \\
&-2nH^2+\frac{n|n-2p|H}{\sqrt{np(n-p)}}\sqrt{
\sum\limits_{\ell=1}^n(h_{\ell\ell}^{(1)})^2-nH^2}\Big)\nonumber \\
&\leq\frac{p(n-p)}{n}\Big(S-2nH^2+\frac{n|n-2p|H}{\sqrt{
np(n-p)}}\sqrt{S-nH^2}\Big) \nonumber\\
&=\frac{p(n-p)}{n}\Big(\|\Phi\|^2+\frac{n|n-2p|}{\sqrt{
np(n-p)}}H\|\Phi\|-nH^2\Big).
\end{align}
From the above inequality we have
\begin{align}\label{in5u}
&\sum_{i=1}^p\sum_{j=p+1}^n\big(2\|\a_f(e_i,e_j)\|^2
-\<\a_f(e_i,e_i),\a_f(e_j,e_j)\>\big)-cp\big(n-p\,{\rm {sgn}}(c)\big)\nonumber\\
&\leq\frac{p(n-p)}{n}\Big(\|\Phi\|^2+\frac{n|n-2p|}{\sqrt{
np(n-p)}}H\|\Phi\|-nH^2-\frac{cn}{n-p}\big(n-p\,{\rm {sgn}}(c)\big)\Big).
\end{align}

It follows from Lemma \ref{lemn} that the inequality \eqref{1} is equivalent 
to the inequality\begin{equation}\label{Su}
\|\Phi\|^2\leq cn\big(2-{\rm {sgn}}(c)\big)+nH^2-\frac{n(n-2k)}{\sqrt{nk(n-k)}}H\|\Phi\|.
\end{equation}
Moreover, equality holds in \eqref{1} at $x$ if and only if \eqref{Su} holds as equality at $x$.
Hence, from \eqref{in5u} and \eqref{Su} we obtain 
\begin{align}\label{in55u}
&\sum_{i=1}^p\sum_{j=p+1}^n\big(2\|\a_f(e_i,e_j)\|^2-\<\a_f(e_i,e_i),\a_f(e_j,e_j)\>\big)-cp\big(n-p\,{\rm {sgn}}(c)\big)\nonumber\\
&\leq \frac{p(n-p)}{n}\Big(cn\big(1-{\rm {sgn}}(c)\big)\frac{n-2p}{n-p}+\sqrt{n}H\|\Phi\|\big(\frac{|n-2p|}{\sqrt{p(n-p)}}-\frac{n-2k}{\sqrt{k(n-k)}}\big)\Big).
\end{align}

If $k=n/2$, then $p=n/2$ and \eqref{in55u} immediately yields \eqref{LS}. 
If $k<n/2$, then \eqref{in55u} is equivalently written as
\begin{align}\label{in6u}
&\sum_{i=1}^p\sum_{j=p+1}^n\big(2\|\a_f(e_i,e_j)\|^2-\<\a_f(e_i,e_i),\a_f(e_j,e_j)\>\big)-cp\big(n-p{\rm {sgn}}(c)\big)\nonumber\\ 
&\leq cp(n-2p)\big(1-{\rm {sgn}}(c)\big)\nonumber\\ 
&+\frac{n^2H\|\Phi\|(k-p)(n-k-p)\sqrt{p(n-p)}}{\sqrt{nk(n-k)}\left(|n-2p|\sqrt{k(n-k)}+(n-2k)\sqrt{p(n-p)}\right)}\nonumber\\ 
&\leq cp(n-2p)\big(1-{\rm {sgn}}(c)\big)\leq 0. 
\end{align}
This completes the proof of the inequality \eqref{LS} for any integer 
$p$ with $k\leq p\leq n/2$ and $p>1$.

Clearly, if the inequality \eqref{1} becomes strict at $x$, so does inequality 
\eqref{Su}. Then \eqref{in55u}, and consequently \eqref{in6u}, become 
strict inequalities. Therefore, \eqref{LS} holds as strict inequality at $x$ 
for any integer $k\leq p\leq n/2$ with $p>1$ and any orthonormal basis 
of $T_xM$ if so does \eqref{1} at $x$.
\vspace{1ex}

Hereafter, we suppose that equality holds in \eqref{LS} 
for a certain orthonormal basis $\{e_i\}_{1\leq i\leq n}$ of $T_xM$ and an integer 
$k\leq p\leq n/2$ with $p>1$. By the above argument we know that 
equality holds in \eqref{1} at $x$. 

We distinguish two cases according to the dimension of the manifold.
\vspace{1ex}

\noindent\emph{Case $n\geq5$}. Since inequality holds in \eqref{LS}, 
all inequalities from \eqref{in1u} 
to \eqref{in6u} become equalities. In particular, from 
 \eqref{in1u}, \eqref{in2u} and Lemma \ref{lem} we have 
\be\label{hu}
h^{(\a)}_{11}=\dots=h^{(\a)}_{pp}
=\lambda_\a,\;\; h^{(\a)}_{p+1p+1}=\dots=h^{(\a)}_{nn}
=\mu_\a \;\,\text{for all}\;\,1\leq\a\leq m.
\ee
Moreover, \eqref{in3u} and \eqref{in4u} imply that 
\begin{align}
h^{(\a)}_{i\tilde i}&=0\;\;\text{for all}\;\;1\leq i\neq \tilde i\leq p, \label{hhu1}\\
h^{(\a)}_{j\tilde j}&=0\;\;\text{for all}\;\;p+1\leq j\neq\tilde j\leq n,\label{hhu2}\\
h^{(\a)}_{ij}&=0\;\;\text{for all}\;\;1\leq i\leq p,\;p+1\leq j\leq n \label{hhu3}
\end{align}
and for any $1\leq\a\leq m$. Observe that the last of the above equalities holds since 
$p(n-p)>n$ for any $2\leq p\leq n-2$. In addition, we have 
\be\label{rsu}
(n-2p)Hh^{(\a)}_{rs}=0\;\,\text{for all}\;\,1\leq r,s\leq n\;\;\text{and any}\;\;2\leq\a\leq m.
\ee
Finally, if $k<n/2$ then \eqref{in6u} yields

\be\label{cp}
cp(n-2p)\big(1-{\rm {sgn}}(c)\big)=0
\ee
and
\be\label{HSu}
H\|\Phi\|(k-p)(n-k-p)=0.
\ee

It follows from \eqref{hhu1}-\eqref{hhu3} that the basis $\{e_i\}_{1\leq i\leq n}$ 
of $T_xM$ diagonalizes the second fundamental form and thus the normal 
bundle of $f$ is flat at $x$. Moreover, using \eqref{hu} we obtain 
\begin{align*}
\a_f(X,T)&=\<X,T\>\eta_1\;\;{\text {for all}}\;\;T\in{\rm {span}}\{ e_1,\ldots,e_p\}, \\
\a_f(X,T)&=\<X,T\>\eta_2\;\;{\text {for all}}\;\;T\in{\rm {span}}\{ e_{p+1},\ldots,e_n\},
\end{align*}
for any $X\in T_xM$, where $\eta_1=\sum_\a\lambda_\a\xi_\a$ and
$\eta_2=\sum_\a\mu_\a\xi_\a$.

Hence, if $f$ is not umbilical at $x$, then both vectors $\eta_1$ and $\eta_2$ 
are distinct Dupin principal normals and this proves part $(i1)$. 

If $f$ is umbilical at $x$, and since equality holds in \eqref{1} at $x$, 
Lemma \ref{lemn} implies that either $c=0$ and $f$ is totally geodesic at 
$x$, or $c=-1$ and $H(x)=\sqrt{3}$. This proves part $(i2)$.

In the sequel, we assume that $k<n/2$. If $c=-1$, then it follows from 
\eqref{cp} and \eqref{HSu} that $n=2p$ and $(H\Phi)(x)=0$. Clearly, 
$f$ is umbilical at $x$ since by assumption $H(x)\geq\sqrt{3}$. That 
$H(x)=\sqrt{3}$ follows from Lemma \ref{lemn} and the fact that equality
 holds in \eqref{1} at $x$. This proves part $(i3)$. 

Now suppose that $c\geq0$. First assume that $(H\Phi)(x)\neq0$. 
Then \eqref{HSu} yields $p=k$ since $k<n/2$. Moreover, it follows from 
\eqref{rsu} that $A_\a=0$ for any $2\leq\a\leq m$. Hence, 
$N_1f(x)={\rm {span}}\left\{\xi_1\right\}={\rm {span}}\left\{\mathcal H(x)\right\}$, 
and this proves part $(i4)$.

Next assume that $(H\Phi)(x)=0$. 
If $\Phi(x)=0$, then from part $(i2)$ we obtain $c=0$ and $f$ is totally 
geodesic at $x$. Suppose now that $\Phi(x)\neq0$, and thus $H(x)=0$. 
By part $(i1)$, there are two distinct Dupin principal normals at $x$ and 
consequently $1\leq\dim N_1f(x)\leq2$. That $H(x)=0$ implies that 
$\dim N_1f(x)=1$ and thus there is a unit vector $\xi\in N_fM(x)$ such that 
$N_1f(x)={\rm {span}}\left\{\xi\right\}$. Moreover, $A_\xi$ has two distinct 
eigenvalues $\lambda$ and $\mu$ with multiplicities $p$ and $n-p$, 
respectively. We may assume that $\lambda\geq0$. Since equality holds 
in \eqref{1} at $x$, we have $S(x)=cn$ and thus $c=1$. Then, a direct 
computation shows that $\lambda=\sqrt{(n-p)/p}$ and 
$\mu=-\sqrt{p/(n-p)}$, and this completes the proof of part $(i5)$.
\vspace{1ex}

\noindent\emph{Case $n=4$}. Since inequality holds in \eqref{LS}, 
all inequalities from \eqref{in1u} to \eqref{in6u} become equalities. 
Then, from \eqref{in1u}, \eqref{in2u}, \eqref{in3u}, \eqref{in4u} and 
Lemma \ref{lem}, we have 
$$
h^{(\a)}_{11}=h^{(\a)}_{22}
=\lambda_\a,\; h^{(\a)}_{33}=h^{(\a)}_{44}
=\mu_\a,\; h^{(\a)}_{12}=h^{(\a)}_{34}=0 \;\,\text{for all}\;\,1\leq\a\leq m.
$$
Now part $(ii)$ easily follows with $\eta_1=\sum_\a\lambda_\a\xi_\a$ and 
$\eta_2=\sum_\a\mu_\a\xi_\a$.\qed

\section{Proofs of the main results}

We recall some facts about Codazzi tensors. A self adjoint endomorphism of 
the tangent bundle of a Riemannian manifold $M$, equipped 
with the Levi-Civit\'a connection $\n$, is said to be a Codazzi tensor if 
$(\n_XA)Y=(\n_YA)X$ for all $X,Y\in\mathcal X(M)$. A smooth distribution 
$E$ on $M$ is \emph{totally geodesic} if $\n_TS\in\Gamma(E)$ whenever 
$T,S\in\Gamma(E)$. The following is well known (see for instance \cite{Der}). 

\begin{lemma}\label{Codaz}
Let $\lambda$ be an eigenvalue of a Codazzi tensor $A$ with constant multiplicity. 
The following assertions hold:\item[(i)] Both $\lambda$ and the corresponding 
eigenbundle $E_\lambda=\{X\in TM: AX=\lambda X\}$ are smooth. Moreover, 
the eigenbundle $E_\lambda$ is integrable and its leaves are totally umbilical.
\item[(ii)] If the multiplicity of the eigenvalue $\lambda$ is greater than one, then 
$\lambda$ is constant along the leaves of the eigenbundle $E_\lambda$.
\item[(iii)] If the eigenvalue $\lambda$ is constant, then the eigenbundle 
$E_\lambda$ is totally geodesic.
\end{lemma}

The following lemma, which will be used repeatedly in the sequel, is an immediate consequence of Proposition 1 in \cite{Reck}.

\begin{lemma}\label{smooth}
Let $f\colon M^n\to\Q_c^{n+m}$ be an isometric immersion with flat normal bundle. Suppose that at each point $x\in M^n$, there exist distinct principal normals $\eta_1(x),\eta_2(x) \in N_fM(x)$ of multiplicities $k$ and $n-k$, respectively, where $k\geq1$ is a fixed integer, such that 
$$
T_xM=E_{\eta_1}(x)\oplus E_{\eta_2}(x).
$$ 
Then the maps 
$x\mapsto\eta_i(x)$, for $i=1,2$, define smooth normal vector fields $\eta_i\in\Gamma(N_fM)$, and the maps $x\mapsto E_{\eta_i}(x)$, for $i=1,2$, 
give rise to a smooth distributions $E_{\eta_1}$ and $E_{\eta_2}$ of ranks 
$k$ and $n-k$, respectively. 
\end{lemma}

We need the following auxiliary results.

\begin{proposition}\label{k<}
Let $f\colon M^n\to\Q_c^{n+m}$, 
$n\geq 5, c\geq0$, be a substantial isometric immersion of a compact 
manifold such that the inequality \eqref{1} is satisfied at any point for 
an integer $2\leq k<n/2$. Assume that equality holds in \eqref{LS} 
for $p=k$ and a certain orthonormal basis of the tangent space 
at any point. Then $M^n$ is isometric to a torus 
$\mathbb T^n_k(r)$ with $r\geq\sqrt{k/n}$, and
 $f$ is the standard embedding into $\Sf^{n+1}$.
\end{proposition}

\proof 
It follows from Proposition \ref{propu} that equality holds in 
\eqref{1} at any point. We consider the open subset 
$M_*=\left\{x\in M^n: (H\Phi)(x)\neq0\right\}$ and we 
distinguish two cases. 
\vspace{1ex}

\indent\emph{Case I}. Suppose that $M_*\neq\emptyset$. 
Part $(i)$-$(i1)$ of Proposition \ref{propu} implies that the submanifold 
has flat normal bundle on $M_*$, and that there exist two distinct 
Dupin principal normal vector fields of multiplicities $k$ and $n-k$,  
with corresponding principal distributions $E_1$ and $E_2$ on $M_*$,  
such that 
$$
T_xM=E_1(x)\oplus E_2(x)
$$ 
for every $x\in M_*$. It follows from Lemma \ref{smooth} that the 
Dupin principal normals and the corresponding principal 
distributions are smooth. Since $H\neq0$, part $(i4)$ of 
Proposition \ref{propu} shows that 
$N_1f={\rm {span}}\left\{\mathcal H\right\}$ on $M_*$. Clearly 
the shape operator $A_1$ associated to $\xi_1=\mathcal H/H$ 
has two pointwise distinct eigenvalues $\lambda,\mu$ of 
multiplicities $k$ and $n-k$ respectively. Let $U$ be a 
connected component of $M_*$. We choose local 
orthonormal frames $\{e_i\}_{1\leq i\leq n}$ and 
$\{\xi_\a\}_{1\leq\a\leq m}$ in the tangent and the 
normal bundle, respectively, such that $\xi_1=\mathcal H/H$, 
$E_1=\spa\left\{e_1,\dots,e_k\right\}$ and 
$E_2=\spa\left\{e_{k+1},\dots,e_n\right\}$. Then we have
$$
A_1e_i=\lambda e_i,\;A_1e_j=
\mu e_j \;\;{\text {for any}}\;\; 1\leq i\leq k,\; k+1\leq j\leq n.
$$
Using 
$$
k\lambda+(n-k)\mu=nH\;\;{\text{and}}\;\;k\lambda^2+(n-k)\mu^2=S,
$$
we find that 
$$
\lambda=H\pm\frac{\sqrt{n-k}}{\sqrt{nk}}\|\Phi\|\;\;{\text{and}}\;\;
\mu=H\mp\frac{\sqrt{k}}{\sqrt{n(n-k)}}\|\Phi\|.
$$
The fact that equality holds in \eqref{1} and Lemma \ref{lemn} yield 
$\|\Phi\|=r(n,k,H,c)$, being $r(n,k,H,c)$ the nonnegative root of the 
polynomial given by \eqref{pol}. From the above we then obtain 
\begin{align}
\lambda&=H\pm\frac{1}{2k}\left(\sqrt{n^2H^2+4ck(n-k)}-(n-2k)H\right)\label{lam},\\
\mu&=H\mp\frac{1}{2(n-k)}\left(\sqrt{n^2H^2+4ck(n-k)}-(n-2k)H\right). \label{mu}
\end{align}

Equations \eqref{lam} and \eqref{mu} imply that if $\mu(x)=0$ at a point 
$x\in M_*$, then $c=1$ and $H^2(x)=k/2(n-2k)$, whereas 
if $\lambda(x)=0$ then $c=0$. 

We claim that $c=1$. Suppose to the contrary that $c=0$. Then we 
claim that $\lambda\mu=0$ on $M_*$. Arguing indirectly, we assume 
that $(\lambda\mu)(x_0)\neq0$ at some point $x_0\in M_*$, and 
thus $\lambda\mu\neq0$ on an open neighborhood $U_0\subset U$ of $x_0$. 
Since $A_\a=0$ for any $2\leq\a\leq m$, it follows from the Codazzi equation 
\be\label{Co}
(\n_XA_\eta)Y-(\n_YA_\eta)X=
A_{\nabla_X^\perp\eta}Y-A_{\nabla_Y^\perp\eta}X,
\ee
for $\eta=\xi_\a,2\leq\a\leq m$, that $\xi_1$ is parallel in the normal bundle 
along $U_0$. Hence $A_1$ is a Codazzi tensor having two distinct eigenvalues 
$\lambda$ and $\mu$ of multiplicity $k\geq2$ and $n-k>k$. By using \eqref{lam} 
and \eqref{mu}, part $(ii)$ of Lemma \ref{Codaz} implies that $\lambda$ and 
$\mu$ are constant on $U_0$. Moreover, part $(iii)$ of Lemma \ref{Codaz} 
implies that both principal distributions $E_1$ and $E_2$ are totally geodesic 
on $U_0$. Using this, we find that for any $1\leq i\leq k,k+1\leq j\leq n$, 
the curvature tensor $R$ of $M^n$ satisfies 
\begin{align*}
\<R(e_i,e_j)e_j,e_i\>&=\<\n_{e_i}\n_{e_j}e_j,e_i \>-\<\n_{e_j}\n_{e_i}e_j,e_i \>
-\<\n_{[e_i,e_j]}e_j,e_i \>\\
&=e_i(\<\n_{e_j}e_j,e_i \>)-\<\n_{e_j}e_j,\n_{e_i}e_i\>-e_j(\<\n_{e_i}e_j,e_i \>)\\
&+\<\n_{e_i}e_j,\n_{e_j}e_i\>-\<\n_{\n_{e_i}e_j}e_j,e_i \>+\<\n_{\n_{e_j}e_i}e_j,e_i \>\\
&=\<\n_{e_i}e_j,\n_{e_j}e_i\>-\<\n_{\n_{e_i}e_j}e_j,e_i \>+\<\n_{\n_{e_j}e_i}e_j,e_i \>. 
\end{align*}
Now, expanding in terms of the orthonormal frames $\{e_{i'}\}_{1\leq i'\leq k}\subset E_1$ and 
$\{e_{j'}\}_{k+1\leq j'\leq n}\subset E_2$, and using again that 
both distributions $E_1$ and $E_2$ are totally geodesic, we obtain
\begin{align*}
\<R(e_i,e_j)e_j,e_i\>&=\sum_{i'=1}^k\<\n_{e_i}e_j,e_{i'}\>\<\n_{e_j}e_i,e_{i'}\>+\sum_{j'=k+1}^n\<\n_{e_i}e_j,e_{j'}\>\<\n_{e_j}e_i,e_{j'}\>\\
&-\sum_{i'=1}^k\<{\n_{e_i}e_j},e_{i'}\>\<\n_{e_{i'}}e_j,e_i \>-\sum_{j'=k+1}^n\<\n_{e_i}e_j,e_{j'}\>\<\n_{e_{j'}}e_j,e_i \>\\
&+\sum_{i'=1}^k\<\n_{e_j}e_i,e_{i'}\>\<\n_{e_{i'}}e_j,e_i \>+\sum_{j'=k+1}^n\<\n_{e_j}e_i,e_{j'}\>\<\n_{e_{j'}}e_j,e_i \>\\
&=0.
\end{align*}
On the other hand, it follows from the Gauss equation that 
$\<R(e_i,e_j)e_j,e_i\>=\lambda\mu$. Hence $\lambda\mu=0$ on $U_0$, 
which is a contradiction. This proves the claim that $\lambda\mu=0$ on 
$M_*$ and consequently $\lambda=0<\mu$ on $M_*$. Thus the immersion 
$f$ has index of relative nullity $k$ at any point of $M_*$. It follows from 
parts $(i2)$ and $(i5)$ of Proposition \ref{propu} that $f$ is totally geodesic 
at any point in $M\smallsetminus M_*$. Hence, the relative nullity of $f$ is 
at least $k$ at any point of $M^n$. Since $M^n$ is compact, there exists a 
point $x_0\in M^n$ and a vector $\xi\in N_fM(x_0)$ such that the shape 
operator $A_\xi$ is positive definite, which is a contradiction. This 
completes the proof of the claim that $c=1$.

We consider the subset $M_{**}=\left\{x\in M_*: (\lambda\mu)(x)\neq0\right\}$. 
It follows from the Codazzi equation \eqref{Co} for $\eta=\xi_\a,2\leq\a\leq m$, 
that $\xi_1$ is parallel in the normal bundle along $M_{**}$. Hence 
$A_1$ is a Codazzi tensor having two distinct eigenvalues $\lambda$ and 
$\mu$ of multiplicity $k\geq2$ and $n-k>k$. By using \eqref{lam} and 
\eqref{mu}, part $(ii)$ of Lemma \ref{Codaz} implies that $\lambda$ and 
$\mu$ are constant on $U$. Moreover, part $(iii)$ of Lemma \ref{Codaz} 
implies that both principal distributions $E_1,E_2$ are totally geodesic. As 
above we have $\<R(e_i,e_j)e_j,e_i\>=0$ for any $1\leq i\leq k,k+1\leq j\leq n$. 
It then follows from the Gauss equation that $\lambda\mu+1=0$ on $M_{**}$. 
Moreover, for any connected component $U\subset M_{**}$ the codimension 
of $f|_U$ can be reduced to one and thus $f(U)$ is an open subset of 
a torus $\mathbb T^n_k(r)$ in $\Sf^{n+1}$ with $r>\sqrt{k/n}$. 

Now we claim that $ M_{**}=M_*$. Observe that $\mu=0$ and $H^2=k/2(n-2k)$ 
on $M_*\smallsetminus M_{**}$ and thus $A_1$ has constant eigenvalues on 
$M_*\smallsetminus M_{**}$. At first we prove that the subset 
$M_*\smallsetminus M_{**}$ has empty interior. Suppose to the 
contrary that $\mathrm{int} (M_*\smallsetminus M_{**})\neq\emptyset$. 
Since $A_\a=0$, for any $2\leq\a\leq m$, it follows from the Codazzi equation 
\eqref{Co} for $\eta=\xi_1$ that $A_1$ is a Codazzi tensor with constant 
eigenvalues on $\mathrm{int} (M_*\smallsetminus M_{**})\neq\emptyset$. 
Part $(iii)$ of Lemma \ref{Codaz} implies that both principal distributions 
$E_1$ and $E_2$ are totally geodesic. Then $\<R(e_i,e_j)e_j,e_i\>=0$ 
for any $1\leq i\leq k,k+1\leq j\leq n$ and the Gauss equation yields 
$\lambda\mu+1=0$ on $\mathrm{int} (M_*\smallsetminus M_{**})\neq\emptyset$. 
This is clearly a contradiction, and thus $M_*\smallsetminus M_{**}$ has 
empty interior. Since $\lambda\mu+1=0$ on $M_{**}$ and $\mu=0$ on 
$M_*\smallsetminus M_{**}$, by continuity it follows that $M_*=M_{**}$. 

Using that $c=1$, part $(i2)$ of Proposition \ref{propu} implies that there are 
no points where $f$ is umbilical. Hence $H=0$ on $M^n\smallsetminus M_*$. 
Then, part $(i5)$ of Proposition \ref{propu} implies that at any point 
$x\in M^n\smallsetminus M_*$, there exists a unit normal vector $\xi$ 
such that $N_1f(x)={\rm {span}}\left\{\xi\right\}$ and the eigenvalues of the 
shape operator $A_\xi$ are $\sqrt{(n-k)/k}$ and $-\sqrt{k/(n-k)}$ with 
multiplicities $k$ and $n-k$, respectively. On the other hand, for any 
connected component $U$ of $M_*$ we know that $f(U)$ is an open 
subset of a torus $\mathbb T^n_k(r)$ in $\Sf^{n+1}$ with $r>\sqrt{k/n}$. 
In particular, the mean curvature is a positive constant on each connected 
component of $M_*$. Since $M^n$ is connected and compact, by continuity 
we obtain $M_*=M^n$. Thus $M^n$ is isometric to a torus $\mathbb T^n_k(r)$ 
with $r>\sqrt{k/n}$ and $f$ is the standard embedding into $\Sf^{n+1}$. 
\vspace{1ex}

\indent\emph{Case II}. 
Suppose that $M_*=\emptyset$. We claim that $c=1$. Indeed, if $c=0$, 
then it follows from parts $(i2)$ and $(i5)$ of Proposition \ref{propu} that 
$f$ is totally geodesic at any point, which is a contradiction. Thus $c=1$. 
Part $(i2)$ of Proposition \ref{propu} implies that there are no points where 
$f$ is umbilical. Hence part $(i5)$ of Proposition \ref{propu} implies that 
there exists a unit normal vector field $\xi$ such that 
$N_1f={\rm {span}}\left\{\xi\right\}$ and the eigenvalues of $A_\xi$ are 
$\sqrt{(n-k)/k}$ and $-\sqrt{k/(n-k)}$ with multiplicities $k$ and $n-k$, 
respectively. Then, the Codazzi equation \eqref{Co} for any normal vector 
field $\eta$ perpendicular to $\xi$ implies that $\xi$ is parallel in the normal 
bundle and the codimension of $f$ can be reduced to one. 
Consequently (cf. \cite{L}), $M^n$ is isometric to a Clifford torus 
$\mathbb T^n_{k}(\sqrt{k/n})$ and $f$ is the standard minimal 
embedding into $\Sf^{n+1}$.\qed

\begin{proposition}\label{M2}
Let $f\colon M^n\to\Q_c^{n+m}$, 
$n\geq 4$, be a substantial isometric immersion of an 
even dimensional manifold with flat normal bundle. Assume that equality 
holds in \eqref{1} for $k=n/2$ at every point, and that there exist two distinct 
Dupin principal normals $\eta_1$ and $\eta_2$, both of multiplicity $k$, such that $$T_xM=E_{\eta_1}(x)\oplus E_{\eta_2}(x)$$ for every point $x\in M^n$ where $f$ is not umbilical. If $c=-1$, suppose further that $H\geq\sqrt{3}$ everywhere. Then, the following assertions hold:
\item[(i)] If $\dim N_1f(x)=2$ at any point $x\in M^n$, then $c\geq0$ and 
$f(M^n)$ is an open subset of a torus $\Sf^k(r)\times \Sf^k(\sqrt{R^2-r^2})$ 
in a sphere $\Sf^{n+1}(R)\subset\Q_c^{n+2}$ of radius $R<1$ if $c=1$. 
\item[(ii)] Suppose that $\dim N_1f(x)=1$ at any point $x\in M^n$. Then, 
the following facts hold:
\begin{itemize}
\item[(a)] If $c=0$, then $f$ is locally a $k$-cylinder in $\R^{n+1}$ over a 
sphere $\Sf^k(r)$, or a $(k-1)$-cylinder over a submanifold $C^{k+1}$ in 
$\R^{k+2}$ which is a cone over a sphere $\Sf^k(r)$.
\item[(b)] If $c=1$, then $f(M^n)$ is an open subset of a torus 
$\mathbb T^n_k(r)$ in $\Sf^{n+1}$.
\item[(c)] If $c=-1$, then $f(M^n)$ is an open subset of a geodesic sphere 
in $\Hy^{n+1}$ with mean curvature $H=\sqrt{3}$. 
\end{itemize}
\end{proposition}

\proof $(i)$ Clearly, there are no points where $f$ is umbilical and 
$H\neq0$ at any point. By our assumption and Lemma \ref{smooth}, 
the two pointwise distinct Dupin principal normal 
vector fields $\eta_1$ and $\eta_2$ are smooth, with corresponding smooth principal distributions $E_{\eta_1}$ and $E_{\eta_2}$, each of rank $k$. We choose a local orthonormal frame $\{\xi_\a\}_{1\leq\a\leq m}$ in the normal 
bundle, with corresponding shape operators $\{A_\a\}_{1\leq\a\leq m}$, 
such that $\{\xi_1,\xi_2\}$ span the plane subbundle $N_1f$ and 
$\xi_1$ is collinear with the mean curvature vector field. Moreover, 
we may choose a local orthonormal frame $\{e_i\}_{1\leq i\leq n}$ 
in the tangent bundle such that $E_{\eta_1}={\rm {span}}\{e_1,\dots,e_k\}$ 
and $E_{\eta_2}={\rm {span}}\{e_{k+1},\dots,e_n\}$. Then we have 
\begin{align}
A_1e_i&=\lambda e_i,\;\;A_1e_j=\mu e_j,\label{A1}\\
A_2e_i&=\rho e_i,\;\;A_2e_j=-\rho e_j,\label{A2}
\end{align}
for any $1\leq i\leq k$ and $k+1\leq j\leq n$, where 
$\lambda=\<\eta_1,\xi_1\>, \mu=\<\eta_2,\xi_1\>$ and 
$\rho=\<\eta_1,\xi_2\>=-\<\eta_2,\xi_2\>$. Clearly, 
$\rho\neq0$ everywhere. Using \eqref{A1} and \eqref{A2}, 
our assumption that equality holds in \eqref{1} on $M^n$ yields 
\be\label{rho}
\lambda\mu+c\big(2-{\rm {sgn}}(c)\big)=\rho^2.
\ee

Since $A_\a=0$ for any $3\leq \a\leq m$, form the Codazzi equation \eqref{Co} 
for $\eta=\xi_\a,3\leq \a\leq m$, $X=e_i$ and 
$Y=e_{\tilde i}, 1\leq i\neq \tilde i\leq k$, we have 
\be\label{Co1}
\lambda\<\nabla_{e_i}^\perp\xi_1,\xi_\a\>+\rho\<\nabla_{e_i}^\perp\xi_2,\xi_\a\>=0\;\;
{\text {for any}}\;\; 1\leq i\leq k.
\ee
Similarly for $X=e_j$ and $Y=e_{\tilde j}, k+1\leq j\neq \tilde j\leq n$, we obtain
\be\label{Co2}
\mu\<\nabla_{e_j}^\perp\xi_1,\xi_\a\>-\rho\<\nabla_{e_j}^\perp\xi_2,\xi_\a\>=0\;\;
{\text {for any}}\;\; k+1\leq j\leq n.
\ee
Finally for $\eta=\xi_\a,3\leq \a\leq m$, $X=e_i, 1\leq i\leq k$, and 
$Y=e_j, k+1\leq j\leq n$, the Codazzi equation \eqref{Co} gives that 
\be\label{Co3}
\mu\<\nabla_{e_i}^\perp\xi_1,\xi_\a\>-\rho\<\nabla_{e_i}^\perp\xi_2,\xi_\a\>=0
\;\;{\text {for any}}\;\; 1\leq i\leq k
\ee
and
\be\label{Co4}
\lambda\<\nabla_{e_j}^\perp\xi_1,\xi_\a\>+\rho\<\nabla_{e_j}^\perp\xi_2,\xi_\a\>=0
\;\;
{\text {for any}}\;\; k+1\leq j\leq n.
\ee
Using that $\rho H\neq0$, we conclude from \eqref{Co1}-\eqref{Co4} that 
both $\xi_1$ and $\xi_2$ are parallel in the normal bundle, and consequently 
$N_1f$ is a parallel plane subbundle of the normal bundle. Hence, $f(M^n)$ 
is contained in a totally geodesic submanifold $\Q^{n+2}_c$ of $\Q^{n+m}_c$ 
and thus $m=2$. 

The shape operator $A_2$ is a Codazzi tensor with two eigenvalues 
$\rho,-\rho\neq0$ of multiplicity $k\geq2$. It follows from part $(ii)$ of 
Lemma \ref{Codaz} that $\rho$ is constant on $M^n$. Moreover, part 
$(iii)$ of Lemma \ref{Codaz} implies that the principal distributions 
$E_{\eta_1}, E_{\eta_2}$ are totally geodesic. Since $A_1$ is also 
a Codazzi tensor, we have from part $(ii)$ of Lemma \ref{Codaz} that
$e_i(\lambda)=0=e_j(\mu)$ for any $1\leq i\leq k$ and any $k+1\leq j\leq n$. 
Using \eqref{A1}, from $(\n_{e_i} A_1)e_j=(\n_{e_j} A_1)e_i$ we obtain 
$$
e_i(\mu)=(\lambda-\mu)\<\n_{e_j}e_i,e_j\>\,\,{\text {and}}\,\,e_j(\lambda)
=(\mu-\lambda)\<\n_{e_i}e_j,e_i\>
$$ 
for any $1\leq i\leq k$ and any $k+1\leq j\leq n$. Since the principal 
distributions $E_{\eta_1}, E_{\eta_2}$ are totally geodesic, we conclude 
that $e_i(\mu)=0=e_j(\lambda)$ for any $1\leq i\leq k$ and any 
$k+1\leq j\leq n$. Therefore, $\lambda,\mu$ are both constant on 
$M^n$. Thus the normal vector field $\xi=2\rho\xi_1+(\mu-\lambda)\xi_2$ 
is parallel in the normal bundle. In addition, it follows form \eqref{A1} 
and \eqref{A2} that $\xi$ is an umbilical direction with corresponding 
shape operator $A_\xi=2\rho HId$. Hence, $f$ is a composition $f=j\circ g$, 
where $g\colon M^n\to\Q^{n+1}_{\tilde c}$ is an isoparametric hypersurface 
with two distinct principal curvatures both of multiplicity $k\geq2$, and 
$j\colon\Q^{n+1}_{\tilde c}\to\Q^{n+2}_c$ is an umbilical inclusion. Since 
$\lambda+\mu=2H$, using \eqref{rho} we obtain 
$$
\|\xi\|=2\sqrt{H^2+c\big(2-{\rm {sgn}}(c)\big)}.
$$ 
The Gauss equation for the umbilical inclusion $j$ gives that the 
curvature of $\Q^{n+1}_{\tilde c}$ is given by 
$$
\tilde c=c+\frac{\rho^2H^2}{H^2+c\big(2-{\rm {sgn}}(c)\big)}>c. 
$$
Using that $f=j\circ g$ and since \eqref{1} holds as equality at 
any point, we find that the squared length $S_g$ of the second 
fundamental form and the mean curvature $H_g$ of $g$ satisfy
\begin{equation}\label{Sg}
S_g=nc\big(1-{\rm {sgn}}(c)\big)+n\tilde c+2nH_g^2.
\end{equation}

We claim that $c\geq0$. Suppose to the contrary that $c=-1$. 
Then $\tilde c>-1$ and thus $\Q^{n+1}_{\tilde c}$ is a geodesic 
sphere, a horosphere, or an equidistant hypersurface of $\Hy^{n+2}$. 
If $\tilde c>0$, then $f(M^n)$ is contained in a torus 
$\Sf^k(r)\times \Sf^k(\sqrt{R^2-r^2})$ in the geodesic 
sphere $\Q^{n+1}_{\tilde c}=\Sf^{n+1}(R)$ with 
$R=1/\sqrt{\tilde c}$. An easy computation shows 
that such a torus satisfies 
\be\label{Sisop}
S_g=n\tilde c+2nH_g^2, 
\ee
which contradicts \eqref{1}. If $\tilde c=0$, then $f(M^n)$ is 
contained in a cylinder $\Sf^k(R)\times\R^k$ in the horosphere 
$\Q^{n+1}_{\tilde c}$. A direct computation shows that such a 
cylinder satisfies \eqref{Sisop} which contradicts \eqref{1}. 
If $-1<\tilde c<0$, then $f(M^n)$ is contained in a cylinder 
$\Sf^k(R)\times\R^k$ in the equidistant hypersurface 
$\Q^{n+1}_{\tilde c}$. Then, it follows from Theorem 1 
in \cite{Ry} that \eqref{Sisop} holds and this again 
contradicts \eqref{1}.

Hence, $c\geq0$ and $f(M^n)$ is contained in a torus 
$\Sf^k(r)\times \Sf^k(\sqrt{R^2-r^2})$ in a sphere 
$\Q^{n+1}_{\tilde c}=\Sf^{n+1}(R)\subset\Q_c^{n+2}$ 
of radius $R=1/\sqrt{\tilde c}$, with $R<1$ if $c=1$. 
\medskip

$(ii)$ We choose a local orthonormal frame 
$\{\xi_\a\}_{1\leq\a\leq m}$ in the normal bundle 
such that $\xi_1$ spans the line subbundle $N_1f$. 
Then $A_\a=0$ for any $2\leq\a\leq m$. Our assumption 
implies that $A_1$ has at most two distinct eigenvalues 
$\lambda, \mu$ at any point. Clearly, $\lambda$ and 
$\mu$ cannot vanish simultaneously. Since by assumption 
equality holds in \eqref{1} we obtain 
\be\label{rho1}
\lambda\mu+c\big(2-{\rm {sgn}}(c)\big)=0.
\ee

We distinguish three cases:
\vspace{1ex}

\noindent\emph{Case $c=0$}. It follows directly from 
\eqref{rho1} that one of $\lambda, \mu$ vanishes and 
thus $f$ has constant index of relative nullity $k$. Without 
loss of generality, we assume that $\lambda>0$. It is known 
that the relative nullity distribution $\mathcal D_f$ is totally geodesic. 
Let $C_T\colon\mathcal D_f^\perp\to\mathcal D_f^\perp$ be the associated 
splitting tensor for any $T\in\mathcal D_f$ (see \cite[p. 186]{DT}). Since 
$$
\a_f(X,Y)=\lambda\<X,Y\>\xi_1\;\;{\text{for all}}\;\;X\in\mathcal X(M^n), Y\in\mathcal D_f^\perp,
$$
from the Codazzi equation
$$
(\n_X^\perp\a_f)(Y,T)=(\n_T^\perp\a_f)(X,Y)\;\;{\text{for all}}\;\;X,Y\in\mathcal D_f^\perp, T\in\mathcal D_f
$$
we find that $C_T=\<\n\log\lambda, T\>Id$, where $Id$ is the identity 
map on the conullity distribution $\mathcal D_f^\perp$. This implies that 
the conullity is umbilical and hence integrable. 

Let $\Sigma^k$ be a leaf of $\mathcal D_f^\perp$. Then, it follows from 
Proposition 7.6 in \cite{DT} that $f$ is locally a generalized cone 
over an isometric immersion $g\colon\Sigma^k\to \Q_{\tilde c}^{k+m}$ 
such that $f\circ j=i\circ g$, where $j\colon\Sigma^k\to M^n$ is the 
inclusion and $i\colon\Q_{\tilde c}^{k+m}\to\R^{n+m}$ is an umbilical inclusion. 
The second fundamental form of the immersion $f\circ j$ is given by 
$$
\a_{f\circ j}(X,Y)=\<X,T\>(\lambda\xi_1+f_*\n\log\lambda)\;\;{\text{for all}}\;\;X,Y\in\mathcal X(\Sigma^k).
$$
Hence, $f\circ j$ is umbilical and $g(\Sigma^k)$ is a sphere $\Sf^k(r)$ 
centered at a point $x_0\in\R^{n+m}$ with radius 
$r=1/\sqrt{\lambda^2+\|\n\log\lambda\|^2}$. If the umbilical 
submanifold $\Q_{\tilde c}^{k+m}$ is totally geodesic in $\R^{n+m}$, 
then the submanifold $f$ is a $k$-cylinder over the sphere $\Sf^k(r)$. 
If $\Q_{\tilde c}^{k+m}$ is a sphere centered at a point 
$\tilde x_0\in\R^{n+m}$, then $f$ is a $(k-1)$-cylinder over a 
submanifold $C^{k+1}$ which is a cone over the sphere 
$\Sf^k(r)$ with vertex at a point $\tilde x_0\neq x_0$.
\vspace{1ex}

\noindent\emph{Case $c=1$}. Condition \eqref{rho1} implies that there are no 
points where $f$ is umbilical. By our assumption and Lemma \ref{smooth}, 
we have two pointwise distinct, smooth Dupin principal normal vector 
fields given by $\eta_1=\lambda\xi_1$ and $\eta_2=\mu\xi_1$, both of multiplicity $k$. 
Hence, $\lambda\neq \mu$ everywhere. We choose a local orthonormal frame 
$\{e_i\}_{1\leq i\leq n}$ in the tangent such that 
\begin{equation*}\label{A3}
A_1e_i=\lambda e_i,\;A_1e_j=\mu e_j \;\;{\text {for any}}\;\; 1\leq i\leq k,\; k+1\leq j\leq n.
\end{equation*}
The Codazzi equation \eqref{Co}
for $\eta=\xi_\a,2\leq \a\leq m$, $X=e_i, 1\leq i\leq k$, and $Y=e_j, k+1\leq j\leq n$, gives that 
$$
\mu\<\nabla_{e_i}^\perp\xi_1,\xi_\a\>=0=\lambda\<\nabla_{e_j}^\perp\xi_1,\xi_\a\>. 
$$
Using 
\eqref{rho1} we conclude that $\xi_1$ is parallel in the normal bundle, and 
consequently $N_1f$ is a parallel line subbundle of the normal bundle. Hence, the 
codimension of $f$ is reduced to one. The hypersurface $f$ has two distinct principal 
curvatures $\lambda,\mu$ both of multiplicity $k\geq2$ satisfying \eqref{rho1}. 
Then, it follows by part $(ii)$ of Lemma \ref{Codaz} and \eqref{rho1} that $\lambda,\mu$ 
are constant, and consequently $f$ is an isoparametric hypersurface. Thus, $f(M^n)$ 
is an open subset of a torus $\mathbb T^n_k(r)$ in $\Sf^{n+1}$. 
\vspace{1ex}

\noindent\emph{Case $c=-1$}. 
We claim that $f$ is umbilical. 
Suppose to the contrary that the open subset 
$M_*=\left\{x\in M^n:\lambda(x)\neq \mu(x)\right\}$ 
is nonempty and let $U$ be a connected component of it. 
Arguing as in Case $c=1$ above, we conclude that the line 
bundle $N_1f$ is a parallel subbundle of the normal bundle 
and thus $f(U)$ is contained as a hypersurface in a totally 
geodesic $\Hy^{n+1}$ of $\Hy^{n+m}$ with two distinct principal 
curvatures $\lambda,\mu$ both of multiplicity $k\geq2$ satisfying 
\eqref{rho1}. Then, part $(ii)$ of Lemma \ref{Codaz} implies that 
$\lambda,\mu$ are constant on $U$, and consequently $f|_U$ 
is an isoparametric hypersurface. It follows from Theorem 1 in \cite{Ry} that 
$\lambda\mu=1$, which contradicts \eqref{rho1} and proves our claim. Since $f$ is 
umbilical, it follows from \eqref{rho1} that $H=\sqrt{3}$, and consequently $f(M^n)$ 
is an open subset of a geodesic sphere of $\Hy^{n+1}$. \qed
\vspace{1.5ex}

\noindent\emph{Proof of Theorem \ref{th1}:}
According to Proposition \ref{propu}, the inequality 
\eqref{LS} is satisfied for any integer $k\leq p\leq n/2$ at any point 
of $M^n$ and for any orthonormal basis of the tangent space at
that point. 

First, suppose that 
\be\label{Hom}
H_p(M^n;\Z)=0=H_{n-p}(M^n;\Z)\;\;{\text {for all}}\;\;k\leq p\leq n/2. 
\ee
Notice that Theorem \ref{er} implies that this is necessarily the case if 
the inequality \eqref{1} is strict at some point of $M^n$.
Since $H_{k-1}(M^n;\Z)$ is a finitely generated abelian group
it decomposes as
$$
H_{k-1}(M^n;\Z)=\Z^{\beta_{k-1}(M)}
\oplus {\rm {Tor}}(H_{k-1}(M^n;\Z)),
$$
where ${\rm {Tor}}(\,)$ denotes the torsion subgroup. By the 
universal coefficient theorem of cohomology, the torsion subgroups of 
$H_{n-k}(M^n;\Z)$ and $H^{n-k+1}(M^n;\Z)$ are isomorphic 
(cf.\ Corollary 4 in \cite[p.\ 244]{Sp}). Since $H_{n-k}(M^n;\Z)=0$ it follows that
${\rm {Tor}}(H^{n-k+1}(M^n;\Z))=0$. Then the Poincar\'e duality
yields ${\rm {Tor}}(H_{k-1}(M^n;\Z))=0$. Thus,
we conclude that $H_{k-1}(M^n;\Z)=\Z^{\beta_{k-1}(M)}$. 

Suppose now that \eqref{Hom} does not hold. Consider the nonempty set
$$
P=\left\{k\leq p\leq n/2: H_p(M^n;\Z)
\neq0\;\;{\text {or}}\;\;H_{n-p}(M^n;\Z)\neq0\right\}
$$
and set $p_*=\max P$. Clearly, $k\leq p_*\leq n/2$ and 
\be\label{p*}
H_{p_*}(M^n;\Z)\neq 0\;\;{\text{or}}\;\; H_{n-p_*}(M^n;\Z)\neq 0.
\ee
Then, it follows from Theorem \ref{er} that at any point $x\in M^n$ 
there exists an orthonormal basis $\{e_i\}_{1\leq i\leq n}$ 
of $T_xM$ such that equality holds in \eqref{LS} for $p=p_*$. Moreover, 
Proposition \ref{propu} implies that \eqref{1} holds as equality 
at any point. 

\indent\emph{Claim}. If $c=-1$, then $k=n/2$. 

Suppose to the contrary that $c=-1$ and $k<n/2$. Then part $(i3)$ of 
Proposition \ref{propu} implies that $f$ is umbilical with $H=\sqrt{3}$. 
It follows from the Gauss equation that $M^n$ has constant sectional 
curvature $K=2$. Hence, $f(M^n)$ is a geodesic sphere in 
$\Hy^{n+1}$ and $M^n$ is isometric to a sphere of radius 
$1/\sqrt{2}$. This contradicts \eqref{p*} thus proving the claim.

We distinguish three cases:
\vspace{1ex}

\noindent\emph{Case I}. Assume that $p_*>k$. From the above 
claim, we have $c\geq0$. Then part $(i4)$ of Proposition \ref{propu} 
implies that $f$ is minimal or umbilical at any point. We claim that $f$ 
is nowhere umbilical. Indeed, it follows from part $(i2)$ of Proposition 
\ref{propu} that points where $f$ is umbilical are totally geodesic points and 
$c=0$. Hence, if there are points where $f$ is umbilical, then $f$ has 
to be a minimal submanifold in the Euclidean space, which is clearly a 
contradiction. Thus there are no points where $f$ is umbilical. Therefore 
$f$ is minimal and $c=1$. Then part $(i5)$ of Proposition \ref{propu} implies 
that there exists a unit normal vector field $\xi$ such that 
$N_1f={\rm {span}}\left\{\xi\right\}$. Moreover, the eigenvalues of $A_\xi$ 
are $\sqrt{(n-p_*)/p_*}$ and $-\sqrt{p_*/(n-p_*)}$ with multiplicities $p_*$ and 
$n-p_*$, respectively. The Codazzi equation \eqref{Co}, if applied to any 
normal vector field $\eta$ perpendicular to $\xi$, implies that $\xi$ is parallel 
in the normal bundle and thus the codimension of $f$ can be reduced to one. 
Consequently (cf. \cite{L}), $M^n$ is isometric to a Clifford torus 
$\mathbb T^n_{p_*}(\sqrt{p_*/n})$ and $f$ is the standard minimal 
embedding into $\Sf^{n+1}$.
\vspace{1ex}

\noindent\emph{Case II}. Suppose that $p_*=k<n/2$. From the above 
claim, we have $c\geq0$. Then, it follows from Proposition \ref{k<} that 
$M^n$ is isometric to a torus $\mathbb T^n_k(r)$ with $r\geq\sqrt{k/n}$ 
and $f$ is the standard embedding into $\Sf^{n+1}$.
\vspace{1ex}

\noindent\emph{Case III}. Suppose now that $p_*=k=n/2$. We know 
from part $(i)$ of Proposition \ref{propu} that $f$ has flat normal 
bundle and that $M^n$ is the disjoint union of the subsets 
$$
M_i=\left\{x\in M^n: \dim N_1f(x)=i\right\},\;\;i=0,1,2.
$$

Clearly, $f$ is nowhere umbilical on the open subset $M_2$. According to 
part $(i1)$ of Proposition \ref{propu} there exist two distinct Dupin principal 
normals of multiplicity $k$ at any point in $M_2$. Moreover, there exist two 
distinct Dupin principal normals of multiplicity $k$ 
at any point in $M_1$ where $f$ is not umbilical. In addition, 
since $M_0$ is nothing but the subset of points where $f$ is totally geodesic, 
it follows directly from part $(i2)$ of Proposition \ref{propu} that if $M_0$ is 
nonempty, then $c=0$. 

We claim that $c\geq0$. Suppose to the contrary that $c=-1$. Part $(i2)$ of 
Proposition \ref{propu} implies that $M_0=\emptyset$. It follows from part 
$(i)$ of Proposition \ref{M2} that $M_2$ is empty and thus $M^n=M_1$. 
Then, part $(ii)$ of Proposition \ref{M2} implies that $f$ is umbilical with 
mean curvature $H=\sqrt{3}$. Hence, $M^n$ is isometric to a sphere of 
radius $1\sqrt{2}$, which contradicts \eqref{p*} and proves the claim.
\vspace{1ex}

\noindent\emph{Case $c=0$}. 
We claim that $M_2$ is nonempty. If otherwise, then part $(ii)$ of 
Proposition \ref{M2} gives that the index of relative nullity of $f$ is 
at least $k$ at any point. This is a contradiction since the 
submanifold is compact and the claim is proved.

It follows from part $(i2)$ of Proposition \ref{propu} that $f$ is nowhere 
umbilical on the open subset $M_1\cup M_2=M^n\smallsetminus M_0$. 
Then, part $(i1)$ of Proposition \ref{propu} together with Lemma \ref{smooth} 
implies that there exist two smooth and distinct Dupin principal normal 
vector fields of multiplicity $k$ on $M_1\cup M_2$. 

Part $(i)$ of Proposition \ref{M2} implies that $S$ is a harmonic function 
on $M_2$. If the interior of $M_1$ is nonempty, then from part $(ii)$ of 
Proposition \ref{M2} we know that the submanifold $f|_{\mathrm{int}(M_1)}$ 
is locally a $k$-cylinder over a sphere $\Sf^k(r)$, and thus $S$ is a 
harmonic function on $\mathrm{int}(M_1)$, or $f|_{\mathrm{int}(M_1)}$ is 
locally a $(k-1)$-cylinder over a cone $C^{k+1}$ shaped over a sphere 
$\Sf^k(r)$. In the later case, $f$ is given locally on $\mathrm{int}(M_1)$ by
$$
f(x,w)=j(x)+w, \;\; (x,w)\in N_i\Sf^{k+m}(R),
$$
where the normal bundle $N_i\Sf^{k+m}(R)$ of the inclusion 
$i\colon\Sf^{k+m}(R)\to\R^{n+m}$ is regarded as subbundle of 
$N_{i\circ j}\Sf^k(r)$ and $j\colon\Sf^k(r)\to\Sf^{k+m}(R)$ is an 
umbilical inclusion with $r<R$. Equivalently, $f$ is locally parametrized by
$$
f(x,t_0,t_1,\dots,t_{k-1})=t_ 0j(x)+\sum_{i=1}^{k-1}t_iv_i,\;t_0>0,\; x\in\Sf^k(r),
$$
where $\Sf^{k+m}(R)\subset\R^{k+m+1}$, $\R^{n+m}=\R^{k+m+1}\oplus\R^{k-1}$ 
and $\{v_i\}_{1\leq i\leq k-1}$ is an orthonormal basis of $\R^{k-1}$. Then, the 
Laplacian operator $\Delta_M$ of $M^n$ is given by
$$
\Delta_M=\frac{k}{t_0}\frac{\partial}{\partial t_0}+\frac{\partial^2}{\partial t_0^2}
+\frac{1}{t^2_0}\,\Delta_{\Sf^k(r)}
+\sum_{i=1}^{k-1}\frac{\partial^2}{\partial t_i^2}.
$$
Since $S=a/t^2_0$, where $a$ is a positive constant, it follows that 
$$\Delta_M S=\frac{a}{t_0^4}(6-n)$$ on ${\rm int}(M_1)$. 
Hence, $S$ is superharmonic (respectively, subharmonic) function on $M^n$ 
if $n\geq6$ (respectively, $n\leq6$), and the maximum principle implies that 
$S$ is a positive constant on $M^n$. Thus $M_0=\emptyset$ and part $(i2)$ 
of Proposition \ref{propu} implies that there are no points where $f$ is umbilical. 

Then, part $(i1)$ of Proposition \ref{propu} and Lemma \ref{smooth} imply 
that there are two pointwise distinct and smooth Dupin principal normal 
vector fields of multiplicity $k$ on $M^n$. Moreover, we claim that one 
of the Dupin principal normal vector fields vanishes on $M_1$. 
Indeed, for any point $x\in M_1$ we choose an orthonormal basis 
$\{\xi_\a\}_{1\leq\a\leq m}$ of $N_fM(x)$ such that 
$N_1f(x)={\rm {span}}\left\{\xi_1\right\}$. Then, the shape 
operators at $x$ satisfy $A_\a=0$ for any $2\leq\a\leq m$ and 
$A_1$ has two distinct eigenvalues $\lambda, \mu$ of multiplicity 
$k$. The fact that equality holds in \eqref{1} yields $\lambda\mu=0$. 
This implies that one of the Dupin principal normal vector fields 
vanishes on $M_1$, whereas the other one is of constant and 
positive length since $S$ is a positive constant. On the other hand, 
it follows from part $(i)$ of Proposition \ref{M2} that both Dupin 
principal normal vector fields have constant and positive length 
on each connected component of $M_2$. By continuity and since 
$M^n$ is compact this implies that $M_2=M^n$. Then, part $(i)$ of 
Proposition \ref{M2} gives that $M^n$ is isometric to a torus 
$\Sf^k(r)\times \Sf^k(\sqrt{R^2-r^2})$, and $f$ is the standard 
embedding into a sphere $\Sf^{n+1}(R)\subset \R^{n+2}$.
\vspace{1ex}

\noindent\emph{Case $c=1$}. Part $(i2)$ of Proposition \ref{propu} 
implies that $M_0=\emptyset$. From part $(i)$ of Proposition \ref{M2} 
we know that if $M_2\neq\emptyset$, then every connected component 
of it is an open subset of a torus in a sphere $\Sf^{n+1}(R)\subset \Sf^{n+m}$ 
of radius $R<1$. Moreover, part $(ii)$ of Proposition \ref{M2} gives that if 
${\rm int}(M_1)\neq\emptyset$, then every connected component of it is 
an open subset of a torus $\mathbb T^n_k(r)$ in a totally geodesic 
$\Sf^{n+1}\subset\Sf^{n+m}$. By continuity, either $M_2=M^n$ or 
$M_1=M^n$, and this completes the proof.\qed
\vspace{1ex}

Let $f\colon M^n\to\Q_c^{n+m}$ be an isometric immersion.
The Gauss equation implies that the Ricci curvature of the 
manifold $M^n$, in the direction of a unit tangent vector 
$X\in T_xM$ at $x\in M^n$, is given by (cf. \cite{DT}) 
\be\label{Ric}
\Ric_M(X)=(n-1)c+\sum_{\a=1}^m(\mathrm{tr}A_\a)
\<A_\a X,X\>-\sum_{\a=1}^m\|A_{\a}X\|^2, 
\ee
where $A_{\a}, 1\leq\a\leq m$, are the corresponding shape 
operators of $f$ associated to an orthonormal 
basis $\{\xi_\a\}_{1\leq\a\leq m}$ of the normal space 
$N_fM(x)$. 

We quote from \cite{OV} the following result which provides 
an estimate for the Bochner operator 
$\B^{[p]}\colon\Omega^p(M^n)\to\Omega^p(M^n)$, a certain 
symmetric endomorphism of the bundle of $p$-forms 
$\Omega^p(M^n)$, in terms of the second fundamental form.

\begin{proposition}\label{lemp}
Let $f\colon M^n\to\Q_c^{n+m}, n\geq 3$, be an isometric immersion. 
The Bochner operator $\B^{[p]}$ of $M^n$, for any $1\leq p\leq n/2$, 
satisfies pointwise the inequality
\be\label{ineqp}
\mathop{\min_{\omega\in \Omega^p(M^n)}}_{\|\omega\|=1}\<\B^{[p]}\omega,\omega\>
\geq \frac{p(n-p)}{n} \big(n(H^2+c) -\frac{n(n-2p)}{\sqrt{np(n-p)}}\, H\|\Phi\|-\|\Phi\|^2\big).
\ee
If equality holds in \eqref{ineqp} at a point $x\in M^n$, then the 
following hold: 
\item[(i)] The shape operator $A_\xi(x)$ has at most two distinct 
eigenvalues with multiplicities $p$ and $n-p$ for every unit vector $\xi\in N_fM(x)$. 
\item[(ii)] If $H(x)\neq 0$ and $p<n/2$, then $N_1f(x)=\mathrm{span}\left\{\mathcal H(x)\right\}$. 
\end{proposition}

Since the Bochner operator $\B^{[1]}$ is nothing but the Ricci tensor,
the inequality \eqref{ineqp} for $p=1$ reduces to an estimate of Ricci 
curvature due to Leung \cite{Le2}. We define the function $\Ric_M^{\min}$ by
$$
\Ric_M^{\min}(x)=\min\left\{\Ric_M(X): X\in T_xM, \|X\|=1\right\},\;x\in M^n.
$$

\begin{proposition}\label{eq}
Let $f\colon M^n\to\Q_c^{n+m}, n\geq3$, be an isometric immersion such 
that the inequality \eqref{1} is satisfied for $k=1$ at any point. If 
$c=-1$, suppose further that the mean curvature satisfies $H\geq\sqrt{3}$ 
everywhere. Then, at any point we have 
\be\label{ineq4}
\Ric_M^{\min}
\geq \frac{n-1}{n} \big(n(H^2+c) -\frac{n(n-2)}{\sqrt{n(n-1)}}\, H\|\Phi\|-\|\Phi\|^2\big)\geq c(n-1)\big(\mathrm{sgn}(c)-1\big).
\ee 
Moreover, if $\Ric_M^{\min}=0$ at any point and $M^n$ is compact, 
then $c=1$, equality holds in \eqref{1} at any point, $M^n$ is isometric to a 
torus $\mathbb T^n_1(r)$ with $r\geq1/\sqrt{n}$ and $f$ is the standard 
embedding into $\Sf^{n+1}$. 
\end{proposition}

\proof
By Lemma \ref{lemn}, our assumption \eqref{1} for $k=1$ is equivalent to 
the inequality
$$
\|\Phi\|^2+\frac{n(n-2)}{\sqrt{n(n-1)}}H\|\Phi\|-n\big(H^2+c\big(2-{\rm {sgn}}(c)\big)\leq0.
$$ 
Then inequality \eqref{ineq4} follows immediately from Proposition \ref{lemp} for $p=1$.

Suppose now that $M^n$ is compact with $\Ric_M^{\min}=0$ at any point. 
Then equality holds in \eqref{ineq4}, and thus $c\geq0$. Moreover, it follows 
from Lemma \ref{lemn} that equality holds in \eqref{1} for $k=1$ at any point. 
Furthermore, Proposition \ref{lemp} implies that the shape operator 
associated to any normal direction has at most two distinct 
eigenvalues with multiplicities $1$ and $n-1$. In addition, we have 
$N_1f(x)=\spa\left\{\mathcal H(x)\right\}$ at points $x\in M_+$, where 
$M_+$ is the subset of points where $H\neq0$.

We claim that $\dim N_1f\leq1$ everywhere on $M^n$. Since the 
shape operator associated to any normal direction has at most two distinct 
eigenvalues, it follows from Lemma 2.2 in \cite{JM} that the eigenvalues of 
each shape operator must have equal multiplicities at points 
$x\in \mathrm{int}(M\smallsetminus M_+)$ with $\dim N_1f(x)>1$. 
This is ruled out for dimension reasons. Thus $\dim N_1f(x)\leq1$ for any 
$x\in \mathrm{int}(M\smallsetminus M_+)$ and consequently by continuity 
$\dim N_1f\leq1$ everywhere on $M^n$. 

Since equality holds in \eqref{1} for $k=1$, the open subset 
$$
M_1=\left\{x\in M^n: \dim N_1f(x)=1\right\}
$$
is nonempty. Let $U_1$ be a connected component of $M_1$.
Then there is a unit normal vector field $\xi$ that spans $ N_1f|_{U_1}$ 
locally. Moreover, we choose a local orthonormal frame 
$\{e_i\}_{1\leq i\leq n}$ in the tangent bundle such that 
$$
A_\xi e_1=\lambda e_1\;\;{\text{and}}\;\; A_\xi e_j=\mu e_j, \;2\leq j\leq n.
$$
Clearly $A_\eta=0$ for any normal vector field perpendicular to $\xi$. 
Thus we have
\begin{align}
\lambda+(n-1)\mu&=nH,\label{HS1}\\
\lambda^2+(n-1)\mu^2&=S=a(n,1,H,c).\label{HS2}
\end{align}
Moreover, it follows from \eqref{Ric} and \eqref{ineq4} that the 
eigenvalues of the Ricci tensor are
\be\label{ric}
\rho_1=c(n-1)+nH\lambda-\lambda^2\geq0
\;\;{\text{and}}\;\;
\rho_2=c(n-1)+nH\mu-\mu^2\geq0, 
\ee
with corresponding multiplicities $1$ and $n-1$, respectively. 
Clearly one of them has to vanish since $\Ric_M^{\min}=0$ everywhere.

We claim that $c=1$. Suppose to the contrary that $c=0$. Then \eqref{HS1} and 
\eqref{HS2} yield 
$$
(\lambda,\mu)=(0,nH/(n-1))\;\;{\text{or}}\;\;(\lambda,\mu)=(2H,(n-2)H/(n-1)) 
$$
at any point of $M^n$. Observe that at points where the latter holds $f$ 
is totally geodesic. Indeed, in this case 
a direct computation shows that
$$
{\rm {Ric}}_M(e_1)=2(n-2)H^2\;\;{\text{and}}\;\;{\rm {Ric}}_M(e_n)
=\frac{n-2}{(n-1)^2}(n^2-2n+2)H^2.
$$
Since one of the eigenvalues of the Ricci tensor has to vanish, 
we conclude that $f$ is totally geodesic at these points. Hence, 
the submanifold has positive index of relative nullity at any point. 
This is contradiction since by the compactness of $M^n$ there is 
a point and a normal direction such that the corresponding shape 
operator is positive definite. 

Thus $c=1$. Observe that there is no totally geodesic point since 
equality holds in \eqref{1} on $M^n$ everywhere for $k=1$. Hence, 
$M_1=M^n$. First claim that the eigenvalues $\rho_1,\rho_2$ of the 
Ricci tensor cannot vanish simultaneously at any point. If otherwise, 
it follows from \eqref{ric} that $\{\lambda, \mu\}=\{\rho_-,\rho_+\}$ at 
a point, where 
$$
\rho_\pm=\frac{1}{2}\big(nH\pm\sqrt{n^2H^2+4(n-1)}\big).
$$
This contradicts \eqref{HS1}. Thus either $\rho_1(x)>0$ or $\rho_2(x)>0$ 
at any point $x\in M^n$. We now claim that $\rho_2>0$ on $M^n$. 
Arguing indirectly suppose that $\rho_2(x_0)=0$ at a point $x_0\in M^n$. 
Then \eqref{ric} imply that 
$\mu(x_0)=\rho_\pm(x_0)$ and $\rho_-(x_0)<\lambda(x_0)<\rho_+(x_0)$. 
If $\mu(x_0)=\rho_-(x_0)$, 
then \eqref{HS1} gives that 
$$
\lambda(x_0)=\frac{1}{2}\big(-n(n-3)H(x_0)+(n-1)\sqrt{n^2H^2(x_0)+4(n-1)}\big),
$$
which contradicts $\rho_-(x_0)<\lambda(x_0)<\rho_+(x_0)$. Similarly, if 
$\mu(x_0)=\rho_+(x_0)$, then \eqref{HS1} gives that 
$$
\lambda(x_0)=-\frac{1}{2}\big(n(n-3)H(x_0)+(n-1)\sqrt{n^2H^2(x_0)+4(n-1)}\big),
$$
which contradicts $\rho_-(x_0)<\lambda(x_0)<\rho_+(x_0)$. Hence, 
$\rho_2>0$ on $M^n$ and consequently $\rho_1=0$ on $M^n$. 

Now we claim that either $\lambda=\rho_-$ at any point of $M^n$ or 
$\lambda=\rho_+$ at any point of $M^n$. Suppose to the contrary 
that there exist points $x_1,x_2$ such that $\lambda(x_1)=\rho_-(x_1)<0$ 
and $\lambda(x_2)=\rho_+(x_2)>0$. By continuity, 
$\lambda$ must vanishes at a point. Since $\rho_1=0$, 
this contradicts the first of \eqref{ric} and proves the claim. 

First suppose that $\lambda=\rho_+$ on $M^n$. Using \eqref{HS1} 
we have 
$$
\mu=\frac{1}{2(n-1)}\big(nH-\sqrt{n^2H^2+4(n-1)}\big)<0.
$$
Then \eqref{HS2} implies that $f$ is minimal. Since $\dim N_1f=1$ 
on $M^n$, arguing as in the proof of Case I in Theorem \ref{th1}, 
we conclude that $M^n$ is isometric to the torus 
$\mathbb T^n_1(1/\sqrt{n})$ and $f$ is the 
standard embedding into $\Sf^{n+1}$.

Now suppose that $\lambda=\rho_-$ on $M^n$. It follows from 
\eqref{HS1} that 
$$
\mu=\frac{1}{2(n-1)}\big(nH+\sqrt{n^2H^2+4(n-1)}\big)>0.
$$
Then, the Codazzi equation \eqref{Co} for any $X=e_1,Y=e_j, 2\leq j\leq n$, 
and any normal vector field $\eta$ perpendicular to $\xi$, implies that $\xi$ 
is parallel in the normal bundle. Hence, the codimension is reduced to 
one and $f$ is contained as a hypersurface in a totally geodesic sphere 
$\Sf^{n+1}$. Observe that $\mu=-1/\lambda$. Then, it follows from 
Theorem 4.2 in \cite{dCD} that the hypersurface is a compact 
rotational hypersurface generated by a unit speed closed curve 
$$
c(s)=\left(x_1(s),0,\dots,0,x_{n+1}(s),x_{n+2}(s)\right) 
$$
with respect to the standard basis $\{\varepsilon_i\}_{1\leq i\leq n+2}$ 
of $\R^{n+2}$. Moreover, the hypersurface is parametrized by
$$
f(t_1,\dots,t_{n-1},s)=\big(x_1(s)\varphi_1,\dots,x_1(s)\varphi_n,x_{n+1}(s),x_{n+2}(s)\big),
$$
where $\varphi(t_1,\dots,t_{n-1})=\left (\varphi_1,\dots, \varphi_n\right)$ 
is a parametrization of the unit sphere in ${\rm {span}} \{\varepsilon_1,\dots, \varepsilon_n\}$. 
Then, the principal curvatures are given by 
$$
\lambda=\frac{\ddot x_1+x_1}{\sqrt{1-x_1^2-\dot x_1^2}},\;\;
\mu=-\frac{\sqrt{1-x_1^2-\dot x_1^2}}{x_1},
$$
with multiplicities $1$ and $n-1$, respectively. Using that $\lambda\mu=-1$, 
we have $\ddot x_1(s)=0$, or equivalently $x_1(s)=as+b,\;a,b\in\R$. Since 
the profile curve is closed, hence periodic, we obtain $a=0$ and thus 
$x_1(s)=b$. Then the principal curvatures are constant. Hence $M^n$ 
is isometric to a torus $\mathbb T^n_1(r)$ with $r\geq1/\sqrt{n}$ and 
$f$ is the standard embedding into $\Sf^{n+1}$.\qed
\vspace{1.5ex}

\noindent\emph{Proof of Theorem \ref{th2}:}
From Lemma \ref{lemn} we know that the inequality \eqref{1} for $k=1$ 
is equivalent to the condition $\|\Phi\|\leq r(n,1,H,c)$, where 
$r(n,k,H,c), 1\leq k\leq n/2$, is the nonnegative root of the 
polynomial given by \eqref{pol}. Since $r(n,1,H,c)\leq r(n,2,H,c)$ 
and $a(n,k,H,c)=(r(n,k,H,c))^2+nH^2$, 
we have 
\begin{equation}\label{SS}
S\leq a(n,1,H,c)\leq a(n,2,H,c).
\end{equation}

We distinguish the following cases according to the dimension.
\vspace{1ex}

\noindent\emph{Case $n\geq 5$}. Theorem \ref{th1} implies 
that either the homology groups satisfy
\begin{equation}\label{H}
H_p(M^n;\Z)=0\;\; {\text {for all}}\;\;2\leq p\leq n-2 \;\; {\text {and}}\;\;H_1(M^n;\Z)=\Z^{\beta_1(M)},
\end{equation}
or equality holds in \eqref{SS} at any point, $c=1$, the manifold 
$M^n$ is isometric to a torus as in parts $(i)$ and $(ii)$ of that 
theorem and $f$ is the standard embedding into $\Sf^{n+1}$.

In the latter case, using Lemma \ref{lemn} 
we have $\|\Phi\|=r(n,1,H,c)=r(n,2,H,c)$ at any point, 
or equivalently 
$$
\|\Phi\|^2+\frac{n(n-2)}{\sqrt{n(n-1)}}H\|\Phi\|-n(H^2+1)=0
$$
and
$$
\|\Phi\|^2+\frac{n(n-4)}{\sqrt{2n(n-2)}}H\|\Phi\|-n(H^2+1)=0. 
$$
The above equations imply that $f$ is minimal. Thus $M^n$ is 
isometric to a Clifford torus $\mathbb T^n_p(\sqrt{p/n}),2\leq p\leq n-2$, 
and $f$ is the standard embedding into $\Sf^{n+1}$. 

Suppose now that the submanifold is not a torus as in parts $(i)$ 
and $(ii)$ of Theorem \ref{th1} and consequently \eqref{H} holds. 
From Proposition \ref{eq} we know that the Ricci curvature is 
nonnegative. Moreover, if $\Ric_M^{\min}=0$ at any point, then 
$M^n$ is isometric to a torus $\mathbb T^n_1(r)$ with 
$r\geq1/\sqrt{n}$ and $f$ is the standard embedding 
into $\Sf^{n+1}$. 

Now assume that there exists a point $x\in M^n$ such that 
$\Ric_M^{\min}(x)>0$. Then, a result due to Aubin \cite{A}, implies 
that the manifold $M^n$ admits a metric of positive Ricci curvature. 
By the Bonnet-Myers theorem the fundamental group $\pi _1(M^n)$ of 
$M^n$ is finite, and consequently $H_1(M^n;\Z)$ is also finite. Since 
$H_1(M^n;\Z)=\Z^{\beta_1(M)}$, we obtain $H_1(M^n;\Z)=0$. Moreover, 
we have $H^n(M^n;\Z)\cong\Z$ for $M^n$ is oriented. The universal 
coefficient theorem of cohomology yields $H_{n-1}(M^n;\Z)$ is torsion 
free. By the Poincar\'e duality $\beta_{n-1}(M)=\beta_1(M)=0$ we 
obtain $H_{n-1}(M^n;\Z)=0$. It follows now from \eqref{H} that 
$M^n$ is a homology sphere. 

Now we consider the universal Riemannian covering 
$\pi\colon\tilde M^n\to M^n$ of $M^n$. Since $\pi _1(M^n)$ 
is finite, $\tilde M^n$ is compact. Obviously the isometric 
immersion $\tilde f=f\circ \pi$ satisfies the inequality 
$\tilde S\leq a(n,1,\tilde H,c)\leq a(n,2,\tilde H,c)$, 
where with tilde we denote the corresponding quantities 
for the immersion $\tilde f$. Since the submanifold $\tilde f$ 
is not a torus as in parts $(i)$ and $(ii)$ of Theorem \ref{th1}, 
then that theorem implies that $H_p(\tilde M^n;\Z)=0$ for all 
$2\leq p\leq n-2$. Arguing as before and we conclude that 
$\tilde M^n$ is a homology sphere with fundamental group 
$\pi _1(\tilde M^n)=0$. Therefore $\tilde M^n$ is a homotopy 
sphere. By the resolution of the generalized Poincar\'{e} 
conjecture due to Smale for $n\geq 5$, we deduce that 
$\tilde M^n$ is homeomorphic to $\Sf^n$. Since the homology 
sphere $M^n$is covered by $\Sf^n$, a result due to Sjerve 
\cite{Sj} implies that $\pi _1(M^n)=0$ and consequently $M^n$ 
is homeomorphic to $\Sf^n$. 
\vspace{1ex}

\noindent\emph{Case $n=4$.} At first we assume that $M^4$ 
is simply connected. The universal coefficient theorem of 
cohomology implies that the torsion subgroups of $H_1(M^4;\Z)$ 
and $H^2(M^4;\Z)$ are isomorphic. Since $M^4$ is simply 
connected, we have $H_1(M^n;\Z)=0$ and thus $H^2(M^4;\Z)$ 
is torsion free. Then, it follows from the Poincar\'e 
duality that also $H_2(M^4;\Z)$ is torsion free. Hence, 
$H_2(M^4;\Z)=\Z^{\beta_2(M)}$.

Suppose that $\beta_2(M)=0$ and thus $H_2(M^4;\Z)=0$. Then,  
the universal coefficient theorem of cohomology implies that 
$H^3(M^4;\Z)$ is torsion free, and by the Poincar\'e duality 
$H_3(M^4;\Z)$ is also torsion free. Thus 
$H_3(M^4;\Z)=\Z^{\beta_3(M)}=\Z^{\beta_1(M)}=0$ and 
consequently $M^4$ is a homology sphere.
By the Hurewicz isomorphism theorem $M^4$ is a homotopy 
sphere. Then, the resolution of the generalized Poincar\'{e} 
conjecture due to Freedman implies that
$M^4$ is homeomorphic to $\Sf^4$. 

Suppose now that $\beta_2(M)\neq0$. Using \eqref{SS} and 
Lemma \ref{lemn}, it then follows from Proposition \ref{lemp} 
that the Bochner operator 
$\B^{[2]}\colon\Omega^2(M^4)\to\Omega^2(M^4)$ 
satisfies for any $\omega\in\Omega^2(M^4)$ the inequality
\be\label{B2}
\<\B^{[2]}\omega,\omega\>
\geq(4c+8H^2-S)\|\omega\|^2\geq4c({\rm {sgn}}(c)-1).
\ee
By the Bochner-Weitzenb\"ock formula the Laplacian of any 
$\omega\in\Omega^2(M^4)$ is given by
$$
\Delta\omega=\n^*\n\omega+\B^{[2]}\omega,
$$
where $\n^*\n$ is the rough Laplacian. From this we obtain 
\be\label{boch.formula}
\<\Delta \omega,\omega\>
=\|\n\omega\|^2+\<\B^{[2]}\omega,\omega\>
+\frac{1}{2}\,\Delta \|\omega\|^2.
\ee 
If $\omega$ is a harmonic $2$-form, then it follows from 
the maximum principle, \eqref{B2} and \eqref{boch.formula}
that $\omega$ has to be parallel. Then, equality holds in \eqref{B2} 
at any point for any harmonic $2$-form. Furthermore, equality holds 
in \eqref{SS} at any point and this implies that the submanifold is 
minimal and therefore $c=1$. Since $\beta_2(M)\neq0$, $M^4$ 
supports a nonzero harmonic $2$-form. Then, if follows from 
Proposition \ref{lemp} that the shape operator $A_\xi(x)$ at any 
$x\in M^4$ and for any $\xi\in N_fM(x)$ has at most two distinct 
eigenvalues with multiplicity $2$. If $m>1$, then Theorem 1.5 in 
\cite{JM} implies that the immersion $f$ is the minimal embedding 
of the complex projective plane $\CP^2_{4/3}$ of constant 
holomorphic curvature $4/3$ into $\Sf^7$ due to Wallach. 
If $m=1$ and since $M^4$ is simply connected, it follows from 
\cite{L} that $M^4$ is isometric to a Clifford torus 
$\mathbb T^4_2(1/\sqrt{2})$, and $f$ is the standard minimal 
embedding into $\Sf^5$. This proves the theorem when $M^4$ 
is simply connected. 

Now suppose that $M^4$ is not simply connected. Proposition \ref{eq} 
implies that the Ricci curvature of $M^4$ is nonnegative. 
Moreover, if $\Ric_M^{\min}=0$ at any point, then $M^4$ is isometric 
to a torus $\mathbb T^4_1(r) $ with $r\geq1/2$ and $f$ is the standard 
embedding into $\Sf^5$. 

Now assume that there exists a point $x\in M^4$ where 
$\Ric_M^{\min}(x)>0$. Then, the result due to Aubin \cite{A}, implies that 
the manifold $M^4$ admits a metric of positive Ricci curvature. By the 
Bonnet-Myers theorem the 
fundamental group $\pi _1(M^4)$ is finite. 
We consider the universal Riemannian covering $\pi\colon\tilde M^4\to M^4$ 
of $M^4$. Since $\pi _1(M^4)$ is finite, $\tilde M^4$ is compact. Obviously 
the isometric immersion $\tilde f=f\circ \pi$ satisfies the inequality 
$\tilde S\leq a(4,1,\tilde H,c)\leq a(4,2,\tilde H,c)$, where with tilde 
we denote the corresponding quantities for the immersion $\tilde f$. 
We claim that $\tilde M^4$ is homeomorphic to $\Sf^4$. If otherwise, then 
from the above argument we know that 
either $\tilde M^4$ is isometric to the Clifford torus $\mathbb T^4_2(1/\sqrt{2})$ 
and $\tilde f$ is the standard minimal embedding into $\Sf^5$, or 
$\tilde M^4$ is isometric to the complex projective plane $\CP^2_{4/3}$ of 
constant holomorphic curvature $4/3$ and $\tilde f$ is the Wallach 
embedding into $\Sf^7$. In either case, the covering map 
$\pi\colon\tilde M^4\to M^4$ has to be a diffeomorphism. 
This contradicts the assumption that $M^4$ is not 
simply connected and proves the claim.
\vspace{1ex}

\noindent\emph{Case $n=3$.} Proposition \ref{eq} implies 
that either $M^3$ is isometric to a torus $\mathbb T^3_1(r)$ 
with $r\geq1/\sqrt{3}$ and $f$ is the standard embedding 
into $\Sf^4$, or there exists a point $x\in M^3$ where 
$\Ric_M^{\min}(x)>0$. In the latter case, the result due 
to Aubin \cite{A} implies that $M^3$ admits a metric 
of positive Ricci curvature. It follows from \cite{H} that
 $M^3$ is diffeomorphic to a spherical space form.\qed 
\vspace{1.5ex}

\noindent\emph{Proof of Corollary \ref{k=2}:} By Theorem \ref{th1}, we 
know that either the submanifold is a torus as in parts $(i)$ and $(ii)$ in 
the statement of the corollary, or the homology groups of $M^n$ satisfy 
$$
H_p(M^n;\Z)=0\;\,\text{for all}\;\,2\leq p\leq n-2
\;\,\text{and}\;\, H_1(M^n;\Z)=\Z^{\beta_1(M)}. 
$$
In the later case, and since $\pi_1(M^n)$ is finite, we can argue as in the 
proof of Theorem \ref{th2} to conclude that $M^n$ is homeomorphic to 
$\Sf^n$.\qed

\noindent Theodoros Vlachos\\
University of Ioannina \\
Department of Mathematics\\
45110 Ioannina -- Greece\\
e-mail: tvlachos@uoi.gr

\end{document}